\documentclass[12pt,a4paper]{article}
\usepackage[utf8]{inputenc} 
\usepackage[T1]{fontenc}    
\usepackage{lmodern}        
\usepackage{amsfonts, amssymb, amsmath, amsthm}
\usepackage{color}
\usepackage{graphicx}
\usepackage[lf]{Baskervaldx} 
\usepackage[bigdelims,vvarbb]{newtxmath} 
\usepackage[cal=boondoxo]{mathalfa} 

\usepackage[width=16.00cm, height=24.00cm, left=2.50cm]{geometry}

\allowdisplaybreaks

\title{The curious world of S.O.T.}

\author{Jorge Antezana\thanks{Supported by PIP11220210100954CO, PID2024-160033NB-I00 funded by MICIU/AEIMCIN/AEI/10.13039/501100011033 and by /10.13039/501100011033FEDER, UE.} \and Daniel Carando\thanks{Supported by CONICET PIP 11220200102366 and UBACYT 20020220300242BA} \and Tom\'as Fern\'andez Vidal\thanks{Supported by CONICET PIP 11220200102366}  \and Melisa Scotti\thanks{Supported by UBACYT 20020220300242BA and CONICET PIP 11220200102366}}

\date{}

\usepackage{graphicx,upref,calligra,float}
\usepackage[rflt]{floatflt}
\usepackage[svgnames]{xcolor}
\usepackage{mathrsfs}

\definecolor{labelkey}{RGB}{40,79,230}
\usepackage[colorlinks=true,linkcolor=azul,citecolor=azul]{hyperref}

\newtheorem{fed}{Definition}[section]
\theoremstyle{plain}
\newtheorem{teo}[fed]{Theorem}
\newtheorem*{teo*}{Theorem}
\newtheorem{lem}[fed]{Lemma}
\newtheorem{cor}[fed]{Corollary}
\newtheorem{pro}[fed]{Proposition}

\theoremstyle{definition}
\newtheorem{rem}[fed]{Remark}

\definecolor{azul}{rgb}{0.1,0.6,0.86}
\definecolor{naranja}{RGB}{249,153,96}

\def\bdem{\begin{proof}}
\def\edem{\end{proof}}
\def\ds{\displaystyle}

\def\bm{\left(\begin{array}}

\def\ben{\begin{enumerate}}
\def\een{\end{enumerate}}
\def\barr{\begin{array}}
\def\earr{\end{array}}
 \def\bit{\begin{itemize}}
\def\eit{\end{itemize}}
\def\beq{\begin{equation}}
\def\eeq{\end{equation}}
\def\bdes{\begin{description}}
\def\edes{\end{description}}

\def\fii{\varphi }

\def\w{\omega}
\def\W{\Omega}
\def\z{\zeta}

\def\N{\mathbb{N}}

\def\Q{\mathbb{Q}}

\def\D{\mathbb{D}}
\def\T{\mathbb{T}}

\newcommand{\peso}[1]{ \quad \mbox{  #1 } \quad }

\newcommand{\pint}[1]{\displaystyle \left \langle\, #1 \, \right\rangle}

\newcommand{\hil}{\mathcal{H}}

\newcommand{\mat}{\mathcal{M}_n(\mathbb{C})}

\newcommand{\polyd}{\D^\infty_2}

\newcommand{\dirvecp}[2]{H_{#1}(\polyd,#2)}

\newcommand{\hsot}[1]{H_\infty^{^{sot}}(\mathbb{T}^\infty,B(#1))}

\newcommand{\dvid}[1]{H_\infty(\polyd,#1)}

\begin{document}


\pagestyle{plain}

\title{Characterization of multipliers on vector-valued Hardy spaces}

\maketitle	

\begin{abstract}
This work characterizes the multipliers on vector-valued Hardy spaces over the infinite polydisk and the infinite polytorus, as well as in the context of Dirichlet series. Unlike the scalar-valued setting, where these frameworks are completely analogous reformulations of one another, there are significant  differences in the vector-valued context. We prove that while the space of multipliers on the infinite polydisk is $H_\infty(\mathbb{D}^\infty_2, B(X))$, the situation on the infinite polytorus is distinct; assuming $X$ is separable, the multiplier space can be identified as $H_\infty^{sot}(\mathbb{T}^\infty, B(X))$, consisting of essentially bounded SOT-measurable functions. These spaces coincide  when $X$ possesses the analytic Radon-Nikodym property. Finally, we extend these results to the associated Hardy spaces of Dirichlet series, $\mathcal{H}_p^+(X)$ and $\mathcal{H}_p(X)$, providing characterizations for their respective multiplier spaces.  
\end{abstract}

\footnotetext[0]{\textit{Keywords:} Multipliers, Hardy spaces, Infinite dimensional analysis, Spaces of Dirichlet series\\
\textit{2020 Mathematics subject classification:} Primary: 30H10,42B15,30B50. Secondary: 46G20 }


\section{Introduction}

The main objective of this work is to characterize pointwise multipliers on vector-valued Hardy spaces over the infinite polytorus and the infinite polydisk, as well as in the context of Dirichlet series, and to examine the connection among them. In the scalar-valued context, these settings are completely analogous reformulations of one another; however, as we will see, this is not the case for vector-valued spaces. 
We begin by presenting the necessary  definitions and background. 

Given a set $\Omega$ and a space of scalar-valued functions $\mathcal{F}$ defined on $\Omega$, a function $\varphi:\Omega\to \mathbb{C}$ is \emph{a multiplier of $\mathcal{F}$} if the pointwise product $\varphi f$ belongs to $\mathcal{F}$ whenever $f\in \mathcal{F}$.  Each multiplier $\varphi$ induces naturally a multiplication operator $M_\varphi : \mathcal{F} \to \mathcal{F}$ defined by $M_\varphi(f)= \varphi f$.
The study of (pointwise) multipliers and the operators they induce has been extensively developed in the literature. Prominent examples include  multipliers of $L_p(X,\mu)$ spaces, Bergman spaces, and Hardy spaces \cite{feldman1999pointwise,hartmann2010pointwise,HeLiSe,mccarthy2004hilbert}. We   now recall some known results regarding the latter.

For $1\leq p <\infty$, the Hardy spaces of functions of one variable are $H_p(\mathbb{T})$ and $H_p(\mathbb{D})$. The first one is defined by all the $p$-integrable functions with respect to the Haar measure on $\mathbb{T}=\{w\in \mathbb{C}, \vert w \vert=1\}$, such that the Fourier coefficients 
$$\hat{f}(n) = \int_{\mathbb{T}} f(w)w^{-n} dw=0,$$
for all negative integers $n$. For $\Omega = \mathbb{D}$, the open unit disk, the Hardy space $H_p(\mathbb{D})$ is the Banach space of all holomorphic functions on $\mathbb{D}$ satisfying 
$$\Vert F \Vert_{H_p(\mathbb{D})}^p = \sup\limits_{0<r<1}\int\limits_{\mathbb{T}} \vert F(rw)\vert^p dw <\infty.$$
The spaces of multipliers of $H_p(\mathbb{T})$ and $H_p(\mathbb{D})$ turns out to be  (see, for example,  \cite{vukotic2003analytic}) the spaces 
$$H_\infty(\mathbb{T})=\{f:\mathbb{T}\to \mathbb{C} : f\in L_\infty(\mathbb{T}),\; \hat{f}(n) = 0 \mbox{ for all } n<0\},$$
and $H_\infty(\mathbb{D})$ of bounded holomorphic functions on $\mathbb D$, respectively. There exists a strong connection between the both Hardy spaces, given by the Fatou's theorem. Every function $F\in H_p(\mathbb{D})$ admits a non-tangential limit to the boundary $\mathbb{T}$ at almost every point $w\in \mathbb{T}$. In particular, \begin{equation}\label{eq-limite}
    f(w)=\lim\limits_{r\to 1^-} F(rw)
\end{equation} 
exists for almost all $w\in \mathbb{T}$. The limit function $f$ is measurable on $\mathbb{T}$, considered with the Haar measure, belongs to $H_p(\mathbb{T})$, and for each positive integer $n$ the Fourier coefficient $\hat{f}(n)$ coincides with the $n$-th Taylor coefficient of $F$. 
The radial limit \eqref{eq-limite}  defines an isometric isomorphism from  $H_p(\mathbb{D})$ to $H_p(\mathbb{T})$.

An analogous characterization of multipliers holds for Hardy spaces in infinitely many variables. Let $\mathbb{T}^\infty$ denote the compact group given by the infinite product of the torus, endowed with the product topology, and its normalized Haar measure  $\mu$. We denote by $\mathbb{Z}^{(\mathbb{N})}$ the set of eventually null sequences of integer numbers. Also, for $w = (w_1,w_2,\ldots)$ and $\alpha\in \mathbb{Z}^{(\mathbb{N})}$ we write $$w^\alpha=w_1^{\alpha_1}\cdots w_n^{\alpha_n} \cdots $$ Note that only finitely many factors in the product  differ from $1$. If $f\in L_1(\mathbb{T}^\infty )$, the $\alpha-$th Fourier coefficient is given by 
$$\hat{f}(\alpha)=\int_{\mathbb{T}^\infty} f(w)w^{-\alpha} dw.$$
The Hardy spaces on the infinite polyotorus $\mathbb{T}^\infty$ are defined as 
$$H_p(\mathbb{T}^\infty) = \Big\{ f\in L_p(\mathbb{T}^\infty,\mu) :  \hat{f}(\alpha )=0 \,\mbox{ if } \alpha\not\in \mathbb{N}_0^{(\mathbb{N})}\Big\}.$$

The multiplier space of $H_p(\mathbb{T}^\infty)$ turns out to be $H_\infty(\mathbb{T}^\infty)$. As in the one-variable case, each $f\in H_p(\mathbb{T}^\infty)$ is associated with a power series $\sum_{\alpha\in \mathbb{N}_0^{(\mathbb{N})}} \hat{f}(\alpha) w^\alpha$. This power series converges absolutely on the domain in $\ell_2$ given by
$$\mathbb{D}_2^\infty : = \{z\in \ell_2 : \vert z_n \vert <1 \mbox{ for all } n\in \mathbb{N}\}$$
and, thus, defines a holomorphic function there (see, for example, \cite[Theorem~13.2]{DirSer}). Given $f\in H_p(\mathbb{T}^\infty)$, we can define $F(z) = \sum_{\alpha\in \mathbb{N}_0^{(\mathbb{N})}} \hat{f}(\alpha) z^\alpha$. Then,  $\hat{f}(\alpha)$  coincides with the  $\alpha$-th Taylor coefficient of $F$, and moreover,
\begin{equation}\label{def norma Hp escalar} 
\Vert f \Vert_{H_p(\mathbb{T}^\infty)} = \sup_{n \in \N} \sup_{0<r<1}
\left(\int_{\mathbb{T}^n} \vert F(rw_1, \ldots, rw_n, 0, \ldots) \vert^{p} \ dw_1 \ldots dw_n\right)^{1/p} < \infty.
\end{equation}
The  Hardy space of holomorphic functions on $\mathbb{D}_2^\infty  $ is defined as
$$H_p(\mathbb{D}_2^\infty) :=\{F:\mathbb{D}_2^\infty \to \mathbb{C} : \mbox{ such that the supremum in} \eqref{def norma Hp escalar} \mbox{ is finite}\}.$$
This space is isometrically isomorphic to $H_p(\mathbb{T}^\infty)$ and the isomorphism maps power series to Fourier series with the same coefficients. As a consequence, the space of multipliers of $H_p(\mathbb{D}_2^\infty)$ are given by the space $H_\infty(\mathbb{D}_2^\infty)$, which in turn coincides with the space $H_\infty(B_{c_0})$ (see \cite[Remmark~13.23]{DirSer}). It is worth noting that unrestricted non-tangential limit do not necessarily exist for functions in infinitely many variables. However, in \cite{aleman2019fatou}, Aleman, Olsen and Saksman presented a boundary approach for which the limit exists almost everywhere.

In \cite{HeLiSe}, Hedenmalm, Lindqvist, and Seip, following an idea of Beurling to study the existence of  Riesz  bases of dilations in $L_2(0,1)$, characterized the multipliers of the Hilbert space of Dirichlet series $\mathcal{H}_2$. A Dirichlet series is a formal expression of the form  $$\sum_{n=1}^\infty a_n n^{-s},$$ where the coefficients $a_n$ are complex numbers and $s$ is a complex variable. The Hardy space $\mathcal{H}_2$ consists of those Dirichlet series whose coefficients belong to $\ell_2$, with norm $$\Vert D\Vert_{\mathcal{H}_2} = \Big(\sum_{n=1}^\infty |a_n|^2\Big)^{1/2}.$$

Since every Dirichlet series in $\mathcal{H}_2$ converges in the right half-plane $$\mathbb{C}_{1/2} :=\{s\in \mathbb{C}: \mbox{Re}(s)>1/2\},$$
 and since $D = 1$ belongs to $\mathcal{H}_2$, it follows that every multiplier must itself be a Dirichlet series convergent on $\mathbb{C}_{1/2}.$  The authors in \cite{HeLiSe} proved that if $D$ is a multiplier, then $D$ actually converges in the whole right half-plane $\mathbb{C}_0:=\{s\in \mathbb{C}: \mbox{Re}(s)>0\}$, and defines, there, a bounded holomorphic function. Consequently, the multiplier space of $\mathcal{H}_2$ is
$$\mathcal{H}_\infty :=\Big\{ D(s)=\sum_{n=1}^\infty a_n n^{-s} :\; \Vert D\Vert_{\mathcal{H}_\infty}=\sup\limits_{\text{Re}(s)>0}\vert D(s) \vert <\infty\Big\}.$$

This space of multipliers is isometrically isomorphic to $H_\infty(\mathbb{T}^\infty)$, and therefore to $H_\infty(B_{c_0})$. The isomorphism is given by the Bohr transform (see \eqref{Def Bohr} below). Later, Bayart \cite{bayart2002hardy} defined the Hardy spaces $\mathcal{H}_p $ of Dirichlet series for arbitrary. For each $1\leq p < \infty$, we have
$$\mathcal{H}_p : = \Big\{ D= \sum\limits_{n=1}^\infty a_n n^{-s} : \mbox{ there exists } f\in H_p(\mathbb{T}^\infty) \mbox{ such that } \hat{f}(\alpha) = a_{\mathfrak{p}^\alpha} \mbox{ for all } \alpha\in \mathbb{N}_0^{(\mathbb{N})} \Big\}.$$
The norm is given by
$$\Big\Vert \sum\limits_{n=1}^\infty a_n n^{-s} \Big\Vert_{\mathcal{H}_p} = \Vert f \Vert_{H_p(\mathbb{T}^\infty)},$$
where $f$ satisfies $\hat{f}(\alpha) = a_{\mathfrak{p}^\alpha}$ for every $\alpha\in \mathbb{N}_0^{(\mathbb{N})}$. Thus $\mathcal{H}_p$  isometrically isomorphic to $H_p(\mathbb{T}^\infty)$, and its multipliers are precisely $\mathcal{H}_\infty$.

\subsection*{Description of the main results}

Since our goal is to extend the scalar description of multipliers to Hardy spaces of vector-valued functions, we must first establish a notion of multipliers in this context. Consider a Banach space $X$, and a space  $\mathcal{F}$ of $X$-valued functions defined  on a set $\Omega$. We seek functions $\varphi$  on $\Omega$ that act linearly on $X$, such that the pointwise multiplication preserves the space $\mathcal{F}$. Accordingly, we define a \emph{multiplier} as a function $\varphi : \Omega \to B(X)$ such that for every $f\in \mathcal{F}$, the product defined by $\varphi \cdot f(\w) = \varphi(\w)f(\w)$ also belongs to $\mathcal{F}$.

In Theorem \ref{Descripcion mult hol}, we characterize the multipliers of $H_p(\mathbb{D}^\infty_2,X)$. We find that, analogous to the scalar case,  the space of multipliers is $H_\infty(\mathbb{D}^\infty_2,B(X))$. The situation in the infinite polytorus is different and requires additional  hypothesis, specifically the separability of $X$. For these Banach spaces, the space of multipliers is generally not $H_\infty(\mathbb{T}^\infty,B(X))$, but $H_\infty^{sot}(\mathbb{T}^\infty,B(X))$, the space of essentially bounded SOT measurable functions. This description is provided in Theorem \ref{mult polytoro}. Finally, Theorem  \ref{iso con sot} establishes the connection between the two spaces of multipliers: they are  isometrically isomorphic when $X$ has the Analytic Radon Nikodym property (ARNP).

The study of multipliers on  Hardy spaces of vector-valued Dirichlet series is somewhat more involved. To begin with, there are two spaces of Dirichlet series naturally induced by the Bohr transform: the space $\mathcal{H}_p^+(X)$, arising from $H_p(\mathbb{D}^\infty_2, X)$, and the space $\mathcal{H}_p(X)$, arising from $H_p(\mathbb{T}^\infty, X)$. These spaces are not isometrically isomorphic unless $X$ has the ARNP. Moreover, the half-plane of convergence of the Dirichlet series depends on the Banach space $X$. We give a characterization of the multipliers, for $\mathcal{H}_p^+(X)$ in Theorem \ref{Teo mult ser dir hol} for every Banach space, and for $\mathcal{H}_p(X)$ in Theorem \ref{Teo mult ser dir toro} for a separable Banach space $X$. Then we study the connection between  both spaces of multipliers.

The article is organized as follows: In Section \ref{seccion def}, we provide the formal definitions of vector-valued Hardy spaces, both on the polydisk and on the polytorus, and establish radial limit$ - $type connections between them. In Section \ref{seccion mult func}, we characterize the multipliers of $H_p(\mathbb{D}^\infty_2, X)$ and $H_p(\mathbb{T}^\infty, X)$. In Section \ref{Seccion series de Dir}, we deal with various vector-valued Hardy spaces of Dirichlet series and their multipliers. In the same section, we use Dirichlet series to study the multipliers of formal power series in each of these spaces. This allows us to complete the picture of the multipliers on both the torus and the polydisk.

\medskip

To keep our notation consistent, from now on we will always use capital letters $F$ and $G$ to refer to holomorphic functions and lowercase letters $f$ and $g$ for the case of measurable functions. The letters $D$ and $E$ will be used in the last section for Dirichlet series. 
\section{Hardy spaces}\label{seccion def}

\subsection{Holomorphic functions on the polydisk}

Given the domain $\polyd=\D^\infty\cap \ell_2$, which we regard as an open subset of $\ell_2$, and a Banach space $X$, we say that a function $F : \mathbb{D}^\infty_2 \to X$ is holomorphic if it is Fr\'echet differentiable at every $z\in \mathbb{D}^\infty_2$. In other words, if for each $z\in \mathbb{D}^\infty_2$ there exists a continuous linear mapping $T_z : \mathbb{D}^\infty_2 \to X$ such that
$$\lim\limits_{h\to 0} \frac{F(z+h) - F(z) - T_z(h)}{\Vert h \Vert_{\ell_2}} =0.$$
We denote by $\mathbb{N}_0^{(\mathbb{N})}$ the eventually null sequences on $\mathbb{N}$. Each holomorphic function $F : \mathbb{D}^\infty_2 \to X$ admits a unique family of coefficients $(c_\alpha)_{\alpha\in \mathbb{N}_0^{(\mathbb{N})}}$ in $X$, the Taylor coefficients of $F$, whose $\alpha-$th term  is given by the Bochner integral
\begin{equation}\label{eq Calpha integral}
c_\alpha=\frac{1}{(2\pi i)^{n}}\int_{|z_1|=r_1}\ldots\int _{|z_{n}|=r_n} \frac{F(z_1,\ldots,z_{n},0,0,\ldots)}{z_1^{\alpha_1+1}\ldots z_{n}^{\alpha_n+1}} \,dz_n\ldots dz_1,
\end{equation}
where $n\in\N$ is such that the coordinates of $\alpha$ are zero from position $n$ onward and $0< r_i < 1$ for all $i=1, \ldots, n$. We call a function  $F : \mathbb{D}^\infty_2 \to X$  analytic at $0$ if there exists a family $(b_{\alpha})_{\alpha\in \mathbb{N}_0^{(\mathbb{N})}}$ in $X$ such that $\sum_{\alpha \in \N_{0}^{(\N)}}b_{\alpha} z^{\alpha}$ converges absolutely, and $F(z)$ coincides point-wise with the limit for all $z\in \mathbb{D}^\infty_2.$ 

It is well known that  every  analytic mapping is holomorphic \cite[Theorem 15.57]{DirSer}. Conversely, not every holomorphic function of infinitely many variables is analytic. For example, the inclusion $I:\mathbb{D}_2^\infty \to c_0$ is holomorphic but not analytic, as seen in \cite[Remark 15.58]{DirSer}. 
However, if $F : \mathbb{D}^\infty_2\to X$ is holomorphic and $N\in \mathbb{N}$, then the restriction of $F$ to the first $N$ variables, denoted by $F_N : \mathbb{D}^N \to X$ and given by $F_N(z) = F(z,0)$, is analytic. Even more, if $z\in \mathbb{D}^N$, $(z,0)\in \mathbb{D}^\infty_2$ and $(c_\alpha)_{\alpha\in \mathbb{N}_0^{(\mathbb{N})}}$ are the Taylor coefficients of $F$, then $\sum_{\alpha\in \mathbb{N}_0^{(\mathbb{N})}} c_\alpha z^\alpha$ converges absolutely and
$$F(z,0) = \sum\limits_{\alpha\in \mathbb{N}_0^{(\mathbb{N})}} c_\alpha z^\alpha = \sum\limits_{\alpha\in \mathbb{N}_0^N} c_\alpha z^\alpha.$$
In particular, for $\alpha\in \mathbb{N}_0^N$ the $\alpha-$th Taylor coefficient of $F_N$ coincides with the $\alpha-$th Taylor coefficient of $F$. 

For $1\leq p < \infty,$ we denote by $\dirvecp{p}{X}$ the $X$-valued Hardy space of holomorphic functions defined on $\polyd$ with values in $X$, such that
\begin{equation*}\label{DefHp}
\|F\|_{\dirvecp{p}{X}} = \sup_{n \in \N} \sup_{0<r<1}
\left(\int_{\mathbb{T}^n} \|F(rw_1, \ldots, rw_n, 0, \ldots) \|_{X}^{p} \ dw_1 \ldots dw_n\right)^{1/p} < \infty.
\end{equation*}
If $p=\infty$, the definition is adapted in the standard way.


We denote by $C_\alpha$ the bounded linear operator from $\dirvecp{p}{X}$ into $X$ that assigns to each function $F \in \dirvecp{p}{X} $, with formal series expansion $ F\sim \sum_{\alpha\in \N_0^{(\N)}} c_\alpha \,z^\alpha$, its $\alpha$-th Taylor coefficient
$$
C_\alpha\left(F\right)=c_\alpha.
$$

If $p=2$ and $X= \hil$ is a Hilbert space, then $\dirvecp{2}{\hil}$ turns out to be a Hilbert space with the inner product given by 
\begin{align}
\pint{F,G}_{\dirvecp{2}{\hil}}=\sum_{\alpha\in \N_0^{(\N)}} \pint{C_\alpha(F),C_\alpha(G)}.\label{PI en Hpolyd}
\end{align}
In particular, the following alternative expression is obtained for the norm:
\begin{align}
\|F\|_{\dirvecp{2}{\hil}}^2=\sum_{\alpha\in \N_0^{(\N)}} \|c_\alpha\|^2_{\hil}.
\end{align}

\subsection{Hardy spaces on the polytorus}\label{Seccion polytoro}

From the Fourier point of view, we consider the compact abelian group $\mathbb{T}^\infty$ with the product topology, endowed with its Haar measure $\mu$. Given a Banach space $X$, a function $f: \mathbb{T}^\infty \to X$ is strongly measurable if there exists a sequence $(\varphi_n)_n$ of simple functions such that $f(w) = \lim\limits \varphi_n(w)$ a.e. on $\mathbb{T}^\infty$. On the other hand, if $\gamma \circ f :\mathbb{T}^\infty \to \mathbb{C}$ is measurable for every $\gamma \in X^*$ then we say that $f$ is weakly measurable. By the Gelfand-Pettis theorem, a weakly measurable function $f$ is strongly measurable if and only if it is almost separable valued. That is, there exists a set $\Omega \subset \mathbb{T}^\infty$ such that $f(\Omega)$ is separable and $\mu(\mathbb{T}^\infty \setminus\Omega)= 0$. In the case where $X$ is a separable Banach space, the weak measurability of $f$ implies that $\Vert f(\cdot)\Vert_X : \mathbb{T}^\infty\to \mathbb{C}$ is measurable. We refer the reader to \cite[Chapter~3]{nikolski2002operators} for a full discussion of these results and definitions.

Now, we can define the spaces $L_p(\T^\infty, X)$ for $1\leq p\leq \infty$ spaces
$$
L_p(\mathbb{T}^\infty,X)=\left\{f:\mathbb{T}^\infty\to X \text{ strongly measurable}: \int\limits_{\mathbb{T}^\infty}\|f(\w)\|_X^p \,d\mu<\infty\right\},
$$
and the definition for the case $p=\infty$ is adapted in the standard way. The $\alpha-$th Fourier coefficient for a function $f\in L_1(\mathbb{T}^\infty,X)$ is given by
$$\hat{f}(\alpha) = \int_{\T^\infty} f(w) \ \overline{w^\alpha} \,dw,$$
for every $\alpha\in \mathbb{Z}^{(\mathbb{N})}$, where again the integral is understood in the Bochner sense. For each $1\leq p\leq \infty$, the Hardy space $H_p(\mathbb{T}^\infty,X)$ is the closed subspace of $L_p(\T^\infty, X)$ defined by

\begin{equation*}
H_{p}(\T^{\infty}, X)=\left\{f\in L_p(\T^\infty,X): \ \int_{\T^\infty} f(\w) \ \overline{\w^\alpha} \,d\w =0 \ \ \mbox{for every $\alpha\notin \N_0^{(\N)}$}\right\}.
\end{equation*}

A very important and useful property of the Hardy spaces $H_p(\mathbb{T}^\infty, X)$ is that they can be approximated by much simpler functions: the analytic trigonometric polynomials. These are finite linear combinations of the monomials with coefficients in $X$. More precisely, a \textbf{trigonometric polynomial} with values in a Banach space $X$ is a finite sum of the form
\[
\sum_{\alpha \in \mathbb{\N}^{(\N_0)}} x_\alpha \, w^\alpha, \quad \text{with } x_\alpha \in X.
\]

\begin{pro} \cite[Proposition~24.6]{DirSer} \label{densidad}
The $X$-valued analytic trigonometric polynomials are norm-dense in $H_p(\mathbb{T}^\infty, X)$ for every $1 \leq p < \infty$, and $w\left(L_\infty(\mathbb{T}^\infty, X), L_1(\mathbb{T}^\infty, X^*)\right)$\text{-}dense in $H_\infty(\mathbb{T}^\infty, X)$.
\end{pro}

\medskip

\subsection{The connection between the polydisk and the polytorus}

 As we have already mentioned, for each $1 \leq p \leq \infty$ and every Banach space $X$, each function $F$ in $H_p(\mathbb{D}^\infty_2,X)$ is associated with a formal power series
$$F\sim\sum\limits_{\alpha\in \mathbb{N}_0^{(\mathbb{N})}}C_\alpha(F)z^\alpha.$$ 
In the same way, every $f\in H_p(\mathbb{T}^\infty,X)$ is associated with a formal power series given by
$$f\sim\sum\limits_{\alpha\in \mathbb{N}_0^{(\mathbb{N})}}\hat{f}(\alpha)w^\alpha.$$
This association with formal power series allows us to establish an isometric isomorphism between the two Banach spaces, under certain conditions on the space $X$. More precisely, we say that $X$ has the Analytic Radon-Nik\'odym property (ARNP) if, for every $F\in H_\infty(\mathbb{D},X)$, the limit
$$\lim\limits_{r\to 1^-}F(rw) \text{ exists for almost every } w\in \mathbb{T}.$$
 Based on this definitions the relation between these two Hardy spaces is given by the following theorem (see \cite[Corollary 24.21]{DirSer} ).


\begin{teo}\label{pablobook1}
For every $1 \leq p \leq \infty$, there exists a unique isometric embedding $$\mathrm{P} : H_p(\mathbb{T}^\infty,X) \to H_p(\mathbb{D}^\infty_2,X),$$
such that $c_\alpha(P[f]) = \widehat{f}(\alpha)$ for every  \( f \in H_p(\mathbb{T}^\infty, X) \) and every $\alpha \in \N_0^{(\N)}.$
Moreover, the Banach space $X$ has the ARNP if and only if $P$ is an isometric isomorphism. 
\end{teo}

\begin{rem}\label{B_c0 a D_2}
    Actually, \cite[Corollary 24.21]{DirSer} is stated for $1\le p<+\infty$. Note, however, that for $p=\infty$, the isometric isomorphism of Theorem \ref{pablobook1} actually holds between $H_\infty(\mathbb{T}^\infty,X)$ and $H_\infty(B_{c_0},X)$, with $B_{c_0}$ the unit ball of the space $c_0.$ The last Hardy space is clearly contained on $H_\infty(\mathbb{D}_2^\infty,X)$, this is because the inclusion $\iota: \ell_2 \hookrightarrow c_0$ is holomorphic and if $F\in H_\infty(B_{c_0},X)$ then $F\circ\iota$ is holomorphic, bounded and coincides with $F$ on the dense sets $\mathbb{D}^m$. In fact, 
    \cite[Corollary~6.5]{defantperez_2018} shows that $H_\infty(B_{c_0},X)$ and $H_\infty(\mathbb{D}_2^\infty,X)$ are isometrically isomorphic for every Banach space $X$.
\end{rem}
\medskip


\medskip

Let us now see that, for functions in the image of $\mathrm{P}$, one can define a kind of radial limit to pass from the polydisk to the polytorus. We follow the definitions given in 
\cite{aleman2019fatou} for the scalar-valued case. 

Given a function $f\in H_p(\mathbb{T}^\infty,X)$, with formal series expansion $$f\sim\sum\limits_{\alpha\in \mathbb{N}_0^{(\mathbb{N})}} c_\alpha w^\alpha,$$
we define $\mathbf{r} = (r^n)_n$ for $0<r<1$, and set $$f_\mathbf{r}(w) := \sum_{\alpha\in \N_0^{(\N)}} c_\alpha \,(\mathbf{r}\cdot w)^\alpha=\sum_{\alpha\in \N_0^{(\N)}} (\mathbf{r}^\alpha\,c_\alpha)\,w^\alpha.$$

This is well-defined since the series converges absolutely on $\mathbb{T}^\infty$, due to the estimate  $$\Vert c_\alpha \Vert_X \leq \Vert f \Vert_{H_p(\mathbb{T}^\infty,X)}.$$ In particular, $f_\mathbf{r}$ is bounded and hence $f_\mathbf{r}$ belongs to $H_p(\mathbb{T}^\infty,X)$. Moreover, these radial functions converge to $f$ in $H_p(\mathbb{T}^\infty,X)$ as $r \to 1^-$, as shown in the following lemma (see \cite[Theorem~1]{aleman2019fatou} for the scalar version).

\begin{lem}\label{aproximacion Hp}
Let $X$ be a Banach space, $1 \leq p < \infty$, and $f \in H_p(\mathbb{T}^\infty,X)$. Then $f_{\mathbf{r}}$ converges to $f$ in $H_p(\mathbb{T}^\infty,X)$.
\end{lem}


\bdem
Let $F = \mathrm{P}(f) \in H_p(\mathbb{D}^\infty_2,X)$. Since $F$ is continuous on $\mathbb{D}^\infty_2$,  we have
$$F(\mathbf{r}\cdot w) = \lim\limits_{N\to \infty} F(rw_1,r^2w_2,\ldots,r^Nw_N, 0).$$
Given that $$F(rw_1,r^2w_2, \ldots,r^Nw_N,0)  = \sum\limits_{\alpha\in \mathbb{N}_0^N} (\hat{f}(\alpha)r^\alpha)w^\alpha$$ for every $N\in \mathbb{N},$ we deduce that $F(\mathbf{r}\cdot) = \mathrm{P} (f_\mathbf{r})$. 
Furthermore, 
\begin{align*} 
\Vert f_{\mathbf{r}} \Vert_{H_p(\mathbb{T}^\infty,X)} ^p &= \Vert F(\mathbf{r}\cdot) \Vert_{H_p(\mathbb{D}^\infty_2,X)}^p = \sup\limits_{N\in \mathbb{N}}\sup\limits_{0<\mathbf{s}<1} \int\limits_{\mathbb{T}^N} \Vert F(\mathbf{r}\cdot \mathbf{s} \cdot w) \Vert_X^p \mathrm{d}w\\
&= \sup\limits_{N\in \mathbb{N}}\sup\limits_{0<\mathbf{s}<\mathbf{r}} \int\limits_{\mathbb{T}^N} \Vert F(\mathbf{s} \cdot w) \Vert_X^p \mathrm{d}w \leq \sup\limits_{N\in \mathbb{N}}\sup\limits_{0<\mathbf{s}<1} \int\limits_{\mathbb{T}^N} \Vert F(\mathbf{s} \cdot w) \Vert_X^p \mathrm{d}w \\
&=\Vert F \Vert_{H_p(\mathbb{D}^\infty_2,X)}^p= \Vert f \Vert_{H_p(\mathbb{T}^\infty,X)} ^p.
\end{align*}

If $P$ is an analytic trigonometric polynomial, then $P_{\mathbf{r}}$ converges uniformly to $P$ as $r\to 1^-$, and in particular in $H_p(\mathbb{T}^\infty,X)$. By the density of the analytic trigonometric polynomials in $H_p(\mathbb{T}^\infty, X)$ (Proposition \ref{densidad}) given $\varepsilon>0$ there exists an analytic trigonometric polynomial $P$ such that
\begin{align*} 
\Vert f -f_{\mathbf{r}} \Vert_{H_p(\mathbb{T}^\infty,X)} &\leq \Vert f - P  \Vert_{H_p(\mathbb{T}^\infty,X)} + \Vert P - P_{\mathbf{r}} \Vert_{H_p(\mathbb{T}^\infty,X)} + \Vert P_{\mathbf{r}} - f_{\mathbf{r}} \Vert_{H_p(\mathbb{T}^\infty,X)} \\
&\leq 2 \Vert f-P\Vert_{H_p(\mathbb{T}^\infty,X)} + \Vert P - P_{\mathbf{r}} \Vert_{H_p(\mathbb{T}^\infty,X)}<\varepsilon. 
\end{align*} 
\edem

For the sake of completeness, we see that if the Banach space $X$ has the ARNP, we have almost everywhere convergence of the radial limits. The proof follows the lines of that of \cite[Theorem 1]{aleman2019fatou} and needs some results on vector-valued functions from \cite{DirSer} and \cite{pisier2016martingales}.

{\begin{pro}
    Let $X$ be a Banach space with the ARNP. Then for each $f\in H_p(\mathbb{T}^\infty,X)$ and $\mathrm{P}(f) = F\in H_p(\mathbb{D}^\infty_2,X)$, $$f(w) = \lim\limits_{r\to1^-} F(rw_1,r^2w_2, \ldots)$$
    for almost every $w \in \mathbb{T}^\infty$.
\end{pro}}

{
\bdem 
By  the inclusion of the spaces $H_p(\mathbb{D}^\infty_2,X)$, it is enough to give 
a proof for the case $p=1$. Following the arguments in 
\cite{aleman2019fatou}, and using \cite[Theorem~3.25 and Proposition~3.22]{pisier2016martingales} we can see the following: Let $\mu$ be a vector  measure on $\mathbb{T}^\infty$ taking values in $X$, whose Fourier coefficients satisfy $\hat{\mu}(\alpha) = 0$ for all $\alpha \notin \mathbb{N}_0^{(\mathbb{N})}$, then the associated maximal function 
$$M\mu (e^{i\theta})= \sup\limits_{r\in(0,1)} \Vert\mu_r(e^{i\theta})\Vert_X,$$
is weak-type $(1,1)$. Here, for each $e^{i\theta}\in \mathbb{T}^\infty$, we write $e^{i\theta} = (e^{i\theta_1}, e^{i\theta_2}\ldots)$, and

$$\mu_\xi = \sum\limits_{\alpha\in \mathbb{Z}^{(\mathbb{N})}} \hat{\mu}(\alpha)r^{\vert\alpha^*\vert_1} e^{it(\alpha_1+2\alpha_2+3\alpha_3+\ldots)}e^{i \alpha\cdot\theta},$$
if $\xi = r e^{it} \in \mathbb{D}$, and $\vert\alpha^*\vert_1 = \sum\limits_{n=1}^\infty\vert n\alpha_n\vert$. 

Now we follow the standard argument given in  \cite[Page~578]{DirSer}. We define the function $L:L_1(\mathbb{T}^\infty,X) \times \mathbb{T}^\infty \to \mathbb{C}$  by
$$L(f,w) = \limsup\limits_{r\to 1} \Vert f(w) - f (\mathbf{r}\cdot w ) \Vert_X. $$
If $f \in H_1(\mathbb{T}^\infty,X),$ in particular $f$ belongs to $ L_1(\mathbb{T}^\infty,X)$, and to prove the statement, it is enough to see that $L(f,w)=0$ for almost every $w\in \mathbb{T}^\infty$. Since the result is clear for trigonometric polynomials, then

$$L(f,w) \leq L(P, w) + L(f-P,w) \leq \Vert f(w) -P(w) \Vert_X + \sup\limits_{0<r<1} \Vert f(\mathbf{r}\cdot w ) -P(\mathbf{r}\cdot w )\Vert_X,$$
for all trigonometric polynomial $P$. Given that $\Vert \cdot \Vert_{1,\infty} \leq \Vert \cdot\Vert_1$ and $M$ is weak-type $(1,1)$ we obtain
\begin{align*}
    \Vert L(f,\cdot) \Vert_{1,\infty} &\leq 2 ( \Vert f-P\Vert_{1,\infty} + \Vert \sup\limits_{0<r<1} \Vert f(\mathbf{r}\cdot w ) -P(\mathbf{r}\cdot w  ) \Vert_X\Vert_{1,\infty} \\
    &\leq 2 \Vert f-P \Vert_1 + \widetilde{C}\Vert f-P \Vert_1.
\end{align*}
By the density of the trigonometric polynomials, we have that $\Vert L(f,\cdot) \Vert_{1,\infty} =0$, in particular $L(f,w) = 0$ for almost every $w \in \mathbb{T}^\infty$, as we wanted to see. 
\edem}


\section{Multipliers of vector-valued Hardy spaces}\label{seccion mult func}

\subsection{On the polydisk}


One of the main result form \cite{HeLiSe} can be reformulated as follows: the space of multipliers on the Banach space \( H_2(\mathbb{D}^\infty_2) \) coincides with \( H_\infty(\mathbb{D}^\infty_2) \). 
Here, by a multiplier we mean a function \( F : \mathbb{D}^\infty_2 \to \mathbb{C} \) such that the multiplication operator
\begin{align*}
  T_F \colon H_2(\mathbb{D}^\infty_2) &\to H_2(\mathbb{D}^\infty_2),\\
  G &\mapsto F  G
\end{align*}
is well-defined. This result holds, actually,  for \( H_p(\mathbb{D}^\infty_2) \),  $1\le p\le +\infty$, as can be deduced from the results in
 \cite{bayart2002hardy} for  Hardy spaces of Dirichlet series.

Our goal in this section is to describe the multipliers of the Hardy spaces $H_p(\mathbb{D}^\infty_2,X)$ for a Banach space $X$. In this context, the product between $F: \mathbb{D}^\infty_2 \to B(X)$  and $G: \mathbb{D}^\infty_2 \to X$ is defined as $$F\cdot G(z) = F(z)(G(z)).$$ 

We say that $F: \mathbb{D}^\infty_2 \to B(X)$ is a \textbf{multiplier} of $H_p(\mathbb{D}^\infty_2,X)$ if the associated multiplication operator if well-defined. That is, if $F\cdot G$ belongs to $H_p(\mathbb{D}^\infty_2,X)$ for every $G \in H_p(\mathbb{D}^\infty_2,X). $

The characterization of the multipliers for the mentioned Hardy spaces is given in the next result.

\begin{teo}\label{Descripcion mult hol}
    Given a Banach space $X$ and $1\leq p < \infty$. A function $F: \mathbb{D}^\infty_2 \to B(X)$ defines a multiplication operator $M_F : H_p(\mathbb{D}^\infty_2, X) \to H_p(\mathbb{D}^\infty_2,X)$ if and only if $F\in H_\infty(\mathbb{D}^\infty_2,B(X))$. In that case, we have that $$\Vert M_F \Vert = \Vert F \Vert_{H_\infty(\mathbb{D}^\infty_2,B(X))}.$$
\end{teo}
\bdem
Assume first that $F\in \dvid{B(X)}$ and let $G\in H_p(\mathbb{D}^\infty_2,X)$. Then $z\mapsto F(z)(G(z))$ is holomorphic from $\mathbb{D}^\infty_2$ to $X$ since it can be thought as the composition of the bilinear function $B(X) \times X \to X$ that sends $(T,x)$ to $T(x)$ and the holomorphic function $z \mapsto (F(z), G(z))$. Even more for every $N\in \mathbb{N}$ and $0<r<1$
\begin{align*} 
\int\limits_{r\mathbb{T}^N} \Vert F(rz_1,rz_2,&\ldots,rz_N,0)(G(rz_1,rz_2,\ldots,rz_N,0))\Vert_X^p dz_1\ldots dz_N \\
&\leq \int\limits_{r\mathbb{T}^N}\Vert F(rz_1,rz_2,\ldots,rz_N,0)\Vert_{B(X)}^p \Vert (G(rz_1,rz_2,\ldots,rz_N,0))\Vert_X^p dz_1\ldots dz_N\\
&\leq \Vert F\Vert_{H_\infty(\mathbb{D}^\infty_2,X)}^p\int\limits_{r\mathbb{T}^N} \Vert (G(rz_1,rz_2,\ldots,rz_N,0))\Vert_X^p dz_1\ldots dz_N.
\end{align*}
Then $F\cdot G \in H_p(\mathbb{D}^\infty_2,X)$, $F$ defines a multiplication operator $M_F$ on $H_p(\mathbb{D}^\infty_2,X)$ and
$$
\|M_F(G)\|_{\dirvecp{p}{X}}\leq \|F\|_{\dvid{B(X)}}\ \|G\|_{\dirvecp{p}{X}}.
$$
Let us now suppose that $F: \mathbb{D}^\infty_2 \to B(X)$ defines a multiplication operator on $H_p(\mathbb{D}^\infty_2,X)$, then $F\cdot G \in H_p(\mathbb{D}^\infty_2,X)$ for every $G\in H_p(\mathbb{D}^\infty_2,X)$. Let $x\in X$, and consider the constant function $G_x \equiv x \in H_p(\mathbb{D}^\infty_2,X)$. Then $$M_F(G_x)(\cdot)=F(\cdot)\,G_x(\cdot) = F(\cdot)x \in H_p(\mathbb{D}^\infty_2 , X)$$ which in particular implies that particular  $F(\cdot)$ is holomorphic for every $x \in X$. Therefore, by \cite[Theorem 15.45]{DirSer} $F$ is holomorphic. 

Given $G \in H_p(\mathbb{D}^\infty_2,X)$ and $N\in \mathbb{N}$, the function 
$$z\mapsto F(z_1,\ldots,z_N,0)\,G(z_1,\ldots,z_N,0)$$ 
defined on $\mathbb{D}^N$ belongs to $H_p(\mathbb{D}^N,X)$ and by the vector-valued version of Cole-Gamelin inequality (see \cite{cole1986representing})
\begin{align*} 
\Vert F(z_1,\ldots,z_N,0)\,G(z_1,\ldots,z_N,0) \Vert_X &\leq \Vert F(\cdot,0)\,G(\cdot,0) \Vert_{H_p(\mathbb{D}^N,X)} \prod\limits_{j=1}^N(1-\vert z_j\vert^2)^{-1/p}\\
&\leq \Vert F(\cdot)\,G(\cdot) \Vert_{H_p(\mathbb{D}^\infty_2,X)} \prod\limits_{j=1}^N(1-\vert z_j\vert^2)^{-1/p}\\
&\leq \Vert M_F \Vert\, \Vert G \Vert_{H_p(\mathbb{D}^\infty_2,X)} \prod\limits_{j=1}^N(1-\vert z_j\vert^2)^{-1/p}.
\end{align*}
For each $x\in X$ and $z \in \mathbb{D}^N$ we define $\Phi_z^x(\z) = x\prod\limits_{j=1}^N\frac{(1-\vert z_j\vert^2)^{1/p}}{(1-\overline{z_j}\z_j)^{2/p}}$. This function belongs to $H_p(\mathbb{D}^N,X)$ subspace of $ H_p(\mathbb{D}^\infty_2,X)$ and $$\Vert \Phi_z^x(\z)\Vert_{H_p(\mathbb{D}^N,X)} = \Vert x \Vert_X \left\Vert \prod\limits_{j=1}^N\frac{(1-\vert z_j\vert^2)^{1/p}}{(1-\overline{z_j}\z_j)^{2/p}}\right\Vert_{H_p(\mathbb{D}^N)} = \Vert x \Vert_X.$$ Therefore
\begin{align*}
    \Vert M_F \Vert\, \Vert x \Vert_X \prod\limits_{j=1}^N(1-\vert z_j\vert^2)^{-1/p} &\geq \Vert F(z_1,\ldots,z_N,0)\,\Phi_z^x(z_1,\ldots,z_N) \Vert_X \\
    &= \Vert F(z_1,\ldots,z_N,0)(x) \Vert_X \prod\limits_{j=1}^N(1-\vert z_j\vert^2)^{-1/p}.
\end{align*}
We can then conclude that $\Vert F(z_1,\ldots,z_N,0) \Vert_{B(X)} \leq \Vert M_F \Vert$ and hence $$\Vert F \Vert_{H_\infty(\mathbb{D}^\infty_2,B(X))} = \sup\limits_{N\in \mathbb{N}} \sup\limits_{z\in \mathbb{D}^N} \Vert F(z,0) \Vert_{B(X)} \leq \Vert M_F \Vert.$$
\edem

\begin{rem}\label{Lemita Tecnico 1}
Let $P = \sum a_\alpha z^\alpha$ be a polynomial with $a_\alpha\in X$, and let $ F$ be a function in $H_\infty(\mathbb{D}_2^\infty,B(X))$ with series representation  $$\sum\limits_{\beta} C_\beta z^\beta,$$ where $C_\beta\in B(X)$ are the  monomial coefficients of $F$. Then, restricting to finitely many variables and using the linearity we get that the  $\gamma-$th coefficient of $M_F(P)$ is given by $\sum\limits_{\alpha+\beta = \gamma} C_\beta(a_\alpha)$. In particular, for every $x \in X$, the function $F_x: \polyd \to X$ defined by $F_x(z) = F(z)(x)$ belongs to $\dvid{X}$, and its monomial coefficients are $c_\alpha(x)$. 

Even more, if $X$ is a Banach space with the ARNP and $G \in H_p(\mathbb{D}^\infty_2,X)$, then there exists a sequence of polynomials $(P_N)_N$ that converges to $G$. Therefore, by the continuity of the coefficient operators on $H_p(\mathbb{D}^\infty_2,X)$ and since $C_\beta\in B(X)$, the $\gamma-$th monomial coefficient of $M_F(G)$ is given by $\sum\limits_{\alpha+\beta = \gamma} C_\beta(a_\alpha)$, where $(a_\alpha)_\alpha$ are the monomial coefficients of $G$.
\end{rem}
\subsection{On the polytorus}

 The aim of this section is to characterize the multipliers of the space $H_p(\mathbb{T}^\infty,X)$. Looking at the scalar case and what we have for holomorphic functions, we might expect the space of multipliers to be $H_\infty(\mathbb{T}^\infty,B(X))$. This is indeed the case when $\dim X=n<\infty$ and, therefore,  
 $B(X) \simeq \mat$. Since $X$ and $\mat$ are separable, from Theorem \ref{mult polytoro} below we get that the space of multipliers on $H_p(\mathbb{T}^\infty,X)$ identifies with $H_\infty(\mathbb{T}^\infty,B(X))$, for every $1\leq p <\infty$. 

If $\dim X = \infty$, in general, this is not the case. As an example, let us consider $X$ a Banach space with unconditional basis $\{x_n\}_n$ and dual basis $\{x^*_n\}_n$. We define $f:\mathbb{T}^\infty\to B(X)$ by $f(w) =T_w$, where $T_w(x) = \sum_n w_n x^*_n(x)x_n$. Since $\{x_n\}_n$ is unconditional,  $T_w$ is uniformly bounded and then $$M=\sup\limits_{w\in \mathbb{T}^\infty}\Vert f(w)\Vert_{B(X)}<\infty.$$
If $g\in H_p(\mathbb{T}^\infty,X)$, we can write $g(w) = \sum_n g_n(w) x_n,$ with $g_n(\cdot) = x_n^*(g(\cdot)) \in H_p(\mathbb{T}^\infty)$. Therefore, for $\gamma\in X^*$, the mapping
$$\gamma(f(w)g(w)) = \sum_n w_n g_n(w)\gamma(x_n),$$
is the point-wise limit of measurable functions and, therefore, measurable. Since $X$ is separable,  $f(\cdot)g(\cdot) :\mathbb{T}^\infty\to X$ is separable valued and, hence, strongly measurable. As $f(w)$ is uniformly bounded on $B(X)$, the function $f(\cdot)g(\cdot$ belongs to $L_p(\mathbb{T}^\infty,X)$. Finally, the $\alpha-$th Fourier coefficient of $f(\cdot)g(\cdot)$ is $\sum_n \hat{g_n}(\alpha-e_n)x_n =0$ if $\alpha\not\in \mathbb{N}_0^{(\mathbb{N})}$, so $f(\cdot)g(\cdot) \in H_p(\mathbb{T}^\infty,X)$ and $f$ defines a multiplier.

On the other hand, for $w,\widetilde{w}\in \mathbb{T}^\infty$ we have 
$$\Vert (T_w - T_{\widetilde{w}})e_n \Vert_X = \Vert (w_n-\widetilde{w}_n)e_n\Vert_X = \vert w_n -\widetilde{w}_n \vert \Vert e_n \Vert_X,$$
and, therefore, $$\Vert f(w) - f(\widetilde{w})\Vert_{B(X)} =\Vert T_w -T_{\widetilde{w}}\Vert_{B(X)} \geq \Vert w -\widetilde{w}\Vert_\infty.$$
In particular, if $A\subset \mathbb{T}^\infty$ is a total measure set, then $f(A)$ is not separable. Hence, $f$ is not essentially separable and then it is not a measurable function. We conclude that the space of  multipliers does not coincide with  $H_\infty(\mathbb{T}^\infty,B(X)).$

This example suggests that, for the polytours, the characterization of the multipliers needs some more work.

\bigskip




\medskip




\subsubsection*{The role of the strong operator topology}

One of the main obstacles when considering \( B(X) \)-valued functions  is that the \( B(X) \) is usually non-separable. In particular, strongly measurable and weakly measurable functions are not equivalent. We have seen that $H_p(\mathbb{T}^\infty,X)$ may admit multipliers that are not strongly measurable. For this reason, we will instead rely on a weaker notion of measurability. Recall that SOT and WOT stand, respectively, for the \emph{strong} and \emph{weak} operator topologies.

Recall that SOT and WOT stand, respectively, for the \emph{strong} and \emph{weak operator topologies}.

\medskip

\begin{fed}
Let $(\W,\mathfrak{S})$ be a measure space. We will say that a function $f:\W\to B(X)$ is SOT-measurable (resp. WOT-measurable) if for every $x\in X$ (resp. $x\in X$ and $\gamma\in X^*$) the function 
$w\mapsto f(w)x$ (resp. $w\mapsto \gamma(f(w)x)$) is strongly measurable in the sense described in Section \ref{Seccion polytoro}.
\end{fed}

\medskip

The following results are particular cases of a result by N. Dunford (see \cite[Section~3.11.1,~Theorem~4]{nikolski2002operators}).

\medskip

\begin{teo}
Let $(\W,\mathfrak{S},\mu)$ be a measurable space and $X$ a separable Banach space.
\begin{itemize}
    \item[a)] A function $f:\W\to B(X)$ is SOT-measurable if and only if is WOT-measurable.
    \item[b)] A function $f:\W\to B(X)$ is strongly measurable if and only if is WOT-measurable and there exists a set $\W_0\subseteq \W$ such that $\mu(\W\setminus \W_0)=0$, and $f(\W_0)$ is separable in $B(X)$.
\end{itemize} 
\end{teo}




\medskip



It is clear that if $f:\mathbb{T}^\infty\to B(X)$ is strongly measurable, then $w\to \Vert f(w)\Vert_{B(X)}$ is a measurable scalar function. Therefore, we can consider the essential supremum norm $$\operatorname*{ess\,sup}_{w \in \mathbb{T}^\infty} \|f(w )\|_{B(X)},$$
and, based on this, define the space $H_\infty(\mathbb{T}^\infty,B(X))$. 
But, as we mentioned before, this space is not enough to characterize the multipliers of $H_p(\mathbb{T}^\infty,X)$. We need a space that includes weakly measurable functions. The separability of $X$ and the SOT-measurability of $f$ are sufficient to ensure the measurability of the scalar function $w\to \Vert f(w) \Vert_{B(X)}$ (see \cite[Lemma 3.4]{curto2023circle}). For this reason, from this point on we will assume that $X$ is a separable space and adopt the notion of SOT-measurability. Following \cite{curto2023circle}, we can consider the space of essentially norm bounded functions, but with a weaker notion of the measurability of $f$.



Let us return to our case of interest $\mathbb{T}^\infty$, equipped with the Haar measure and its corresponding $\sigma$-algebra. Then
 $$
L_\infty^{sot}(\mathbb{T}^\infty,B(X))=\left\{f:\T^\infty \to B(X): \operatorname*{ess\,sup}_{w \in \mathbb{T}^\infty} \|f(w)\|_{B(X)} < \infty\ \text{ and } f \text{ is SOT-measurable}\right\},
$$

identifying two functions if they agree almost everywhere with respect to the Haar measure, and the norm of the space given by the essential supremum. 

This space turn out to be a Banach space (see \cite[Lemma~3.2]{curto2023circle}).
Now we can move on to defining a SOT-Hardy space of operator-valued functions on the polytorus.

\begin{fed}
Given a separable Banach space $X$ and $1\leq p\leq \infty$, we define
$$
H_\infty^{sot}(\mathbb{T}^\infty,B(X)):=\{f\in L_\infty^{sot}(\mathbb{T}^\infty,B(X)) : f(\cdot) x \in  H_\infty(\mathbb{T}^\infty,X) \; \mbox{for every } x\in X\},
$$
with the norm given by the essential supremum.
\end{fed}

\medskip

\begin{rem}
In general, the space $H_\infty^{sot}(\mathbb{T}^\infty,B(X))$ is different from $H_\infty(\T^\infty, B(X)),$ 
regardless  of the domain being  $\T$ or  $\T^\infty$ (see  \cite[Example~3.6]{curto2023circle} for $\T$ or the example at the beginning of this section for $\T^\infty$). It is clear that if $B(X)$ is separable then  both spaces coincide, but this hardly happens. A notable exception is  the hereditarily indecomposable space constructed in \cite{Argyros} to solve the scalar-plus-compact problem. For the complexification of this particular space, the equality $H_\infty^{sot}(\mathbb{T}^\infty,B(X)) = H_\infty(\mathbb{T}^\infty,B(X))$ holds.
\end{rem}

As in the case of holomorphic functions, given $f\in\hsot{X}$, then for every $X-$valued function $g$ on $\mathbb{T}^\infty$ we can define
$$
M_f(g)(w)=f(w)(g(w)).
$$
Let us see that these functions are the ones that characterize the multipliers of $H_p(\mathbb{T}^\infty,X)$.


\begin{teo}\label{mult polytoro}
    Given $X$ a separable Banach space and $1\leq p < \infty$. A function $f: \mathbb{T}^\infty \to B(X)$ defines a multiplication operator on $M_f : H_p(\mathbb{T}^\infty, X) \to H_p(\mathbb{T}^\infty,X)$ if and only if $f\in H_\infty^{sot}(\mathbb{T}^\infty,B(X))$. In that case $\Vert M_f \Vert = \Vert f \Vert_{H_\infty^{sot}(\mathbb{T}^\infty,B(X))}$.
\end{teo}


\bdem
Suppose first that $f\in H^{sot}_\infty(\mathbb{T}^\infty,B(X))$. Let us start by seeing that $f\cdot g : \mathbb{T}^\infty \to X$ is strongly measurable for every $g\in H_p(\mathbb{T}^\infty,X)$. Given that $g$ is strongly measurable, there exists a sequence of simple functions $g_n$ such that $g_n\to g$ almost everywhere, and hence $f\cdot g_n \to f\cdot g$ almost everywhere. If $g_n (w)= \sum b_j^n \chi_{A_k}(w)$, then $f(w)(b_j^n)$ is strongly measurable, since $f$ is sot-measurable, and then $$f(w)(g_n(w)) = \sum f(w)(b_j^n) \chi_{A_k}(w),$$
is strongly measurable. Therefore, $f\cdot g$ is the almost everywhere pointwise limit of strongly measurable functions and, then, strongly measurable. The measurability of $w\mapsto\Vert f(w)\Vert_{B(X)}$ and the uniformly boundedness $$\operatorname*{ess\,sup}_{w \in \mathbb{T}^\infty} \|f(w)\|_{B(X)} < +\infty,$$ give  $f\cdot g\in L_{p}(\T^{\infty}, X)$
%
and
$$
\|f\cdot g\|_{L_{p}(\T^{\infty}, X)}\leq \|f\|_{\hsot{X}}\ \|g\|_{L_{p}(\T^{\infty}, X)}.
$$
Hence, $M_f :L_p(\mathbb{T}^\infty,X) \to L_p(\mathbb{T}^\infty,X)$ is well defined and  bounded.
On the other hand, note that by the definition of $\hsot{X}$,
$$
f\cdot\sum_{n=0}^N x_n \w^{\alpha_n} \in H_{\infty}(\mathbb{T}^{\infty}, X) \subset H_p(\mathbb{T}^\infty,X),
$$
for every analytic polynomial. Therefore, by the density of the polynomials in $H_p(\mathbb{T}^\infty,B(X)),$ we see that $M_f$ defines a multiplier on $H_{p}(\T^{\infty}, X)$.

Let us suppose  now that $f: \mathbb{T}^\infty \to B(X)$ defines a multiplication operator on $H_p(\mathbb{T}^\infty,X)$. Then, for each $x\in X$ we can take $g_x \equiv x \in H_p(\mathbb{T}^\infty,X)$ and, therefore, $f(\cdot)x \in H_p(\mathbb{T}^\infty,X)$. Hence, $f$ is a SOT-measurable function.Also, for $\gamma\in X^*$ we have that $\gamma(f(\cdot)x)$ belongs to the scalar Hardy space $H_p(\mathbb{T}^\infty).$ Also, if $h \in H_p(\mathbb{T}^\infty)$, then for every $w\in \mathbb{T}^\infty$ 
$$(\gamma(f(\cdot)x) h)(w) = \gamma(f(w) ( xh(w))),$$
with $xh \in H_p(\mathbb{T}^\infty,X)$. Then
\begin{align*} 
\Vert \gamma(f(\cdot)x) \varphi(\cdot)\Vert_{H_p(\mathbb{T}^\infty)} \leq \Vert \gamma \Vert_{x^*} \Vert f(\cdot) ( xh(\cdot)) \Vert_{H_p(\mathbb{T}^\infty,X)} \leq \Vert \gamma\Vert_{X^*} \Vert M_f \Vert \Vert x \Vert_X \Vert h \Vert_{H_p(\mathbb{T}^\infty)}.
\end{align*}
We see that
$\gamma(f(\cdot)x)$ defines a multiplication operator on $H_p(\mathbb{T}^\infty)$ and, hence, belongs to $H_\infty(\mathbb{T}^\infty)$ and 
$$\Vert \gamma(f(\cdot)x) \Vert_{H_\infty(\mathbb{T}^\infty)} = \Vert M_{\gamma(f(\cdot)x)} \Vert \leq \Vert \gamma\Vert_{X^*} \Vert M_f \Vert \Vert x \Vert_X$$
We then have
\begin{align*}
    \sup\limits_{w\in\mathbb{T}^\infty} \Vert f(w)x\Vert_X & = \sup\limits_{w\in\mathbb{T}^\infty}\sup\limits_{\gamma\in B_{X^*}} \vert \gamma(f(w)x)\vert \leq \Vert M_f \Vert \Vert x \Vert_X,
\end{align*}
from which we conclude  that $f(\cdot)x \in H_\infty(\mathbb{T}^\infty,X)$, and 
\begin{align*}
    \sup\limits_{w\in\mathbb{T}^\infty} \Vert f(w)\Vert_{B(X)} & = \sup\limits_{w\in\mathbb{T}^\infty}\sup\limits_{x\in B_{X^*}} \Vert f(w)x\Vert_X \leq \Vert M_f \Vert.
\end{align*}
We have shown that $f\in H_\infty^{sot}(\mathbb{T}^\infty,B(X))$ and $\Vert f \Vert_{H_\infty^{sot}(\mathbb{T}^\infty,B(X))} \leq \Vert M_f \Vert$.
\edem

\begin{rem}\label{Lemita tecnico 2}
If $X$ is separable and $f \in H_\infty^{sot}(\mathbb{T}^\infty,B(X)),$ then $f$ is formally associated with a Fourier series $\sum a_\alpha w^\alpha$ with $a_\alpha \in B(X)$. The coefficient $a_\alpha$ is given by  $$a_\alpha(x) = \hspace{.2em} \widehat{\hspace{-.2em}f(\cdot)(x)\hspace{-.2em}}\hspace{.2em} (\alpha) = \int_{\T^\infty } f(w)(x) {w}^{-\alpha} dw.$$ The integral makes sense since $ f(\cdot)(x)$ belongs to $H_p(\mathbb{T}^\infty,X)$. Also, it is clear that $$\sup\limits_\alpha \Vert a_\alpha\Vert_{B(X)} \leq \Vert f\Vert_{H^{sot}_\infty(\mathbb{T}^\infty,B(X))}.$$
Moreover, by the linearity of $M_f$ and the density of the polynomials in $H_p(\mathbb{T}^\infty,X)$, we get that if $g\in H_p(\mathbb{T}^\infty,X)$ then the power series of $M_F(g)$ is $\sum c_\alpha z^\alpha$ with $$c_\alpha = \sum\limits_{\beta+\gamma = \alpha} a_\beta(b_\gamma),$$ and $b_\gamma = \hat{g}(\gamma)$.
\end{rem}

In both the case of $H_p(\mathbb{D}_2^\infty,X)$ and $H_p(\mathbb{T}^\infty,X)$ (if $X$ is separable), each multiplier $\varphi$ can be associated, at least formally, with a power series $\sum a_\alpha z^\alpha$, with $a_\alpha\in B(X)$. Working with formal power series, one can define the product of $\sum a_\alpha z^\alpha$ with $a_\alpha\in B(X)$ and $\sum b_\alpha z^\alpha$ with $b_\alpha \in X$ by
$$\left(\sum\limits_{\alpha\in \mathbb{N}_0^{(\mathbb{N})}} a_\alpha z^\alpha\right) \left( \sum\limits_{\alpha\in \mathbb{N}_0^{(\mathbb{N})}} b_\alpha z^\alpha \right) = \sum\limits_{\alpha\in \mathbb{N}_0^{(\mathbb{N})}} \left( \sum\limits_{\beta + \gamma =\alpha}a_\beta(b_\gamma)\right)z^\alpha,$$
being a formal power series with coefficients in $X$. It is a natural question to ask which formal power series  $\sum a_\alpha z^\alpha$ with $a_\alpha\in B(X)$ satisfy that the formal power series 
$$\sum\limits_{\alpha\in \mathbb{N}_0^{(\mathbb{N})}} \left( \sum\limits_{\beta + \gamma =\alpha}a_\beta(b_\gamma)\right)z^\alpha,$$
is associated to a function in $H_p(*,X)$ for each $\sum b_\alpha z^\alpha \sim g\in H_p(*,X)$, being $*=\mathbb{D}_2^\infty$ or $\mathbb{T}^\infty$. We will see in Section \ref{Seccion series de Dir}, working with Dirichlet series, that this power series are characterized by $H_\infty(\mathbb{D}^\infty_2,B(X))$ and $H_\infty^{sot}(\mathbb{T}^\infty,X)$.

\subsection{The connection between multiplier spaces.}
If $X=\mathbb{C}$ and $1\leq p <\infty$, the spaces of multipliers of $H_p(\mathbb{D}^\infty_2,\mathbb{C})$ and $H_p(\mathbb{T}^\infty,\mathbb{C})$ are, respectively, $H_\infty(\mathbb{D}_2^\infty,\mathbb{C})$ and $H_\infty(\mathbb{T}^\infty,\mathbb{C})$. These spaces are isomorphic via the identification of their Cauchy and Fourier coefficients. It is natural to ask whether an analogous isomorphism holds for the multipliers spaces of $X$-valued functions. We will show that this is indeed the case under certain conditions, specifically when the Banach space $X$ is separable and has the ARNP. 
 In order to see this, we define
$$P_\infty: H_\infty^{sot}(\mathbb{T}^\infty,B(X)) \to H_\infty(\mathbb{D}^\infty_2,B(X)),$$
by $F= P_\infty(f)$ if and only if $f$ and $F$ are formally associated to the same power series in the terms of Remarks \ref{Lemita Tecnico 1} and \ref{Lemita tecnico 2}. The good definition of $P_\infty$ is proved in the following result.

\begin{teo}\label{iso con sot}
    Let $X$ be a separable Banach space. Then, $$P_\infty:\hsot{X}\to\dvid{B(X)},$$ is  an into isometry. Moreover, if $X$ has the ARNP then $P_\infty$ is an isomorphism.
\end{teo}

It is worth noticing that the above result is analogous to \cite[Theorem~1.1]{curto2023circle} but for the case infinitely many variables and $p= \infty$.

\medskip


\medskip

\bdem
Consider $f\in\hsot{X}$. Then, for every $x\in X$ we have that  $f(\cdot)x$ belongs to $H_{\infty}(\T^{\infty}, X)$ and hence, there exists $F_x\in \dirvecp{\infty}{X}$ such that
 $$
 C_\alpha(F_x) = \widehat{(f(\cdot)x)}(\alpha)   \peso{and}  \|f(\cdot)x\|_{H_{\infty}(\T^{\infty}, X)}=\|F_x\|_{\dirvecp{\infty}{X}}. 
 $$
 For each $z\in \mathbb{D}_2^\infty$ we define  $F(z):X\to X$ by setting $F(z)(x)= F_x(z)$. By the linearity of $f(w)x$ and of the Fourier transform together with the uniqueness of both monomial and Fourier coefficients, it follows that $F(z)$ is linear. Moreover, 
  \begin{align*}
     \sup_{z\in\mathbb{D}^\infty_2}\sup_{\|x\|_X=1}\| F(z) (x)\|_X &=  \sup_{\|x\|_X=1} \sup_{z\in\mathbb{D}^\infty_2} \| F(z) (x)\|_X =  \sup_{\|x\|_X=1} \sup_{z\in\mathbb{D}^\infty_2} \| F_x(z)\|_X \\ 
     & = \sup\limits_{\Vert x \Vert_X =1} \operatorname*{ess\,sup}_{w \in \mathbb{T}^\infty} \Vert f(w)x\Vert_X =  \operatorname*{ess\,sup}_{w \in \mathbb{T}^\infty} \sup\limits_{\Vert x \Vert_X =1} \Vert f(w)x\Vert_X\\
     &=  \|f\|_{\hsot{X}}.\\
 \end{align*}
This not only shows that the operator $F(z)$ is bounded for every  $z\in \polyd$, but also that $F$ is uniformly bounded on $\polyd$. On the other hand, since for each $x \in X$ the function $z \mapsto F(z)(x)$ is holomorphic, it follows that $F$ is holomorphic. Therefore $P_\infty (f) = F \in H_\infty(\mathbb{D}_2^\infty,B(X))$ and $P_\infty(f) $ is an into isometry.

Let us suppose now that $X$ has the ARNP and let $F\in \dvid{B(X)}$ and $(c_\alpha)\subset B(X)$ with
$$
F(z)\sim\sum_{\alpha\in \N_0^{(\N)}} c_\alpha z^\alpha. 
$$
 First we need to show that there exists a SOT-measurable function $f:\T^\infty \to B(X)$ such that $$
\sum_{\alpha\in \N_0^{(\N)}} c_\alpha w^\alpha \sim f \in \hsot{X},
$$
in the terms of the Remark \ref{Lemita tecnico 2}. 
\medskip


\medskip

With this goal in mind, let us first note that for every $x \in X$ $$
F_x(z)=F(z)x,
$$
is holomorphic and
$$
\|F_x(z)\|_X \leq \|F\|_{\dvid{B(X)}} \|x\|_{X}.
$$
Hence  $F_x\in H_{\infty}(\polyd, X)$. Furthermore, by Remark \ref{Lemita Tecnico 1} we also know that its series expansion is given by
$$
F_x(z)=\sum_{\alpha\in \N_0^{(\N)}} c_\alpha(x) z^\alpha.
$$
By Theorem  \ref{pablobook1}, we deduce that 
$$
f_x(w)=\sum_{\alpha\in \N_0^{(\N)}} c_\alpha(x) w^\alpha
$$
defines a function in $H_{\infty}(\T^\infty, X)$ with the same norm (in other words, there exists a unique element in $H_{\infty}(\T^\infty, X)$ whose Fourier coefficients are the $c_\alpha(x)$, here we use that $X$ has the ARNP). Observe that the Fourier coefficients of $f_x$ are precisely $c_\alpha(x)$. We now aim to combine the information corresponding to each element $x \in X$ preserving the measurability of the function. Let  $X_0= \{x_n\}$ be a dense subset of  $X$. 


Then for each element $x,y\in X_0$ and each coefficient $\alpha,\beta\in\Q$, we can consider the  functions $f_{\alpha x+\beta y}$, $f_x$ and $f_y$. 
Comparing coefficients we get that
$$
f_{\alpha x+\beta y}(w)=\alpha f_x(w) +\beta f_y(w),
$$
for all $w$ in a set $\W(x,y,\alpha,\beta)\subseteq \T^\infty$ of total measure. In this way, it becomes clear that it is possible to obtain a set of full measure $\W \subseteq \T^\infty$ such that for every $x,y\in X_0$,  $\alpha,\beta\in\Q$ and $w\in\W$ we get
\begin{enumerate}
\item[i.)] $f_{\alpha x+\beta y}(\w)=\alpha f_x(\w) +\beta f_y(\w)$. 
\item[ii.)] $\ds \sup_{\w\in\W}\Vert f_x(\w)\Vert_X= \|f_x\|_{H_{\infty}(\T^\infty, X)}$.
\end{enumerate}
Then, for every $x,y\in X_0$ the following inequality holds:
\begin{align}
\|f_x(\w)-f_y(\w)\|_{X}&=\|f_{x-y}(\w)\|_{X}\leq \|f_{x-y}\|_{H_{\infty}(\T^\infty, X)}\nonumber \\&=\|F_{x-y}\|_{H_{\infty}(\polyd, X)}\nonumber\\&\leq \|F\|_{\dvid{B(X)}} \|x-y\|_{X}.\label{eq sale uniforme}
\end{align}
For each $w\in\W$ we define $f(w)x=f_x(w)$ if $x\in X_0$ and if $y\notin X_0$ we choose a sequence $\{y_n\}$ in $X_0$ that converges to $y$ to define 
$$
f(w)y=\lim_{n\to\infty} f_{y_n}(w),
$$
Note that the limit exists and the convergence is uniform on $\W$, due to \eqref{eq sale uniforme}. Moreover, from \eqref{eq sale uniforme}, we can deduce that $f(w)y$ is well-defined, because it is independent of the choice of the sequence $\{y_n\}$, and a bounded linear operator since
$$
\|f(w)y\|_{X}\leq \|F\|_{\dvid{B(X)}} \|y\|_{X}.
$$
Now, if we extend it as zero on the complement of $\W$, then for each $x \in X$ the function defined on  $\T^\infty$ by $w\mapsto f(w)x$ belongs to $H_{\infty}(\T^\infty, X)$. Indeed, this is clear by construction if $x \in X_0$.  Moreover, since the subspace $H_{\infty}(\mathbb{T}^\infty, X)$ is closed, the same conclusion holds for every $y \in X$. Therefore, we conclude that $f$ is SOT-measurable and that $f \in \hsot{X}$. To complete the argument we need to show that 
$$
\widehat{f}(\alpha)=c_\alpha,
$$
in the SOT sense of Remark \ref{Lemita tecnico 2}. From the definition of the Fourier coefficients  it is clear that for every $x\in X_0$ we have $\widehat{f}(\alpha)(x)=c_\alpha(x)$, and by density it holds for every $x \in X$. We conclude that $ P_\infty(f)=F$ and therefore $P_\infty$ is an isometric isomorphism.
\edem

\medskip



From the first part of Theorem \ref{iso con sot} we get the next relation between the multiplication operators.

\begin{teo}\label{dos multiplicadores}
Let $X$ be a separable Banach space,  $f\in\hsot{X}$ and  $F=P_\infty(f)$. Then
\begin{equation}\label{eq compatibildiad}
\mathrm{P} M_f= M_F \mathrm{P},
\end{equation}
as an operator on $H_p(\mathbb{T}^\infty,X)$ for every $1\leq p <\infty$, where  $P$ is the isometric embedding defined in Theorem \ref{pablobook1}.
\end{teo}


\bdem 
It is enough to prove this identity on polynomials, which form a dense subset of the space $H_p(\mathbb{T}^\infty,X)$ for every $1\leq p <\infty$. Moreover, by linearity, it suffices to check it on functions of the form 
$$
\{xw^\alpha: \, x\in X,\, \alpha\in \N_0^{(\N)}\}.
$$

By definition $\mathrm{P}(xw^\alpha)(\z)=x\z^\alpha$ for every $\z\in \mathbb{D}_2^\infty$. Then, if 
$$
f(w)\sim\sum_{\alpha\in \N_0^{(\N)}} c_\alpha w^\alpha,
$$
understanding this expression in the SOT sense, we have
$$
\mathrm{P}_\infty(f)=F\sim\sum_{\alpha\in \N_0^{(\N)}} c_\alpha z^\alpha.
$$
Thus,
$$
M_f(xw^\alpha)\sim\sum_{\tilde \alpha\in \N_0^{(\N)}}c_{\tilde\alpha}(x) w^{\tilde\alpha + \alpha} \peso{and} M_F(\mathrm{P}(xz^\alpha))=\sum_{\tilde \alpha\in \N_0^{(\N)}}c_{\tilde\alpha}(x) z^{\tilde\alpha + \alpha}.
$$
Concluding that \eqref{eq compatibildiad}  is true.
\edem

Combining this result with Lemma \ref{aproximacion Hp} we obtain the following corollary.

\medskip

\begin{cor}\label{Aprox segunda parte}
Given $f\in\hsot{X}$ and  $\fii\in H_p(\mathbb{T}^\infty,X)$, let $F=\mathrm{P}_\infty(f)$ and $\Phi=\mathrm{P}(\varphi)$. For each $r\in (0,1)$, let $\mathbf{r}=\{r^n\}_{n\in\N}$. Then
$$
\lim_{r\to 1^{-}} \int_{\T^\infty} \big\|f(w)(\fii(w))-F(\mathbf{r}\cdot w)(\Phi(\mathbf{r}\cdot w)) \big\|^p_X\ dw = 0,
$$
where $\mathbf{r}\cdot w=(rw_1,r^2w_2,r^3w_3,\ldots)$.
\end{cor}

As a consequence of Theorems \ref{iso con sot} and \ref{pablobook1} we get

\begin{cor}
    Let $X$ a Banach space. If $B(X)$ is separable, then $B(X)$ has the ARNP if and only if $X$ has the ARNP.
\end{cor}

\medskip

\section{Multipliers in vector-valued Dirichlet series}\label{Seccion series de Dir}

In this section, we shift the study of multipliers for vector-valued functions to the perspective of Dirichlet series. Let us start with some basic notions of them. For a Banach space $X$, an $X$-valued Dirichlet series is a formal expression of the form
$$D= \sum\limits_{n=1}^\infty a_n n^{-s},$$
where $a_n\in X$ and $s$ is a complex variable. It is a well known result that if a Dirichlet series $D$ converges on some complex $s_0 = \sigma_0 + it_0$ then the series converges and defines a holomorphic function on the half-plane 
$$\mathbb{C}_{\sigma_0} := \{s\in \mathbb{C} : \text{Re}(s) > \sigma_0\}.$$
In this way, each Dirichlet series determines the abscissas of convergence and absolute convergence defined by
$$\sigma_c(D) := \inf\{\sigma\in \mathbb{R} : D(s) \mbox{ converges on } \mathbb{C}_\sigma\},$$
$$\sigma_a(D) := \inf\{\sigma\in \mathbb{R}: D \mbox{ converges absolutely on } \mathbb{C}_\sigma\},$$
being $\sigma_a(D)\leq \sigma_c(D)+1.$

We start by defining the  formal multiplication for $X$-valued Dirichlet series.

\begin{fed}\label{Def prod dirichlet}
Let $X$ be a Banach space and let $D = \sum\limits_{n} a_n n^{-s}$ and $E= \sum\limits_n b_n n^{-s}$ be Dirichlet series with coefficients in $B(X)$ and $X$ respectively. We define $D\cdot E$ as the Dirichlet series with coefficients in $X$ given by
$$D\cdot E = \sum\limits_n \left(\sum\limits_{k\cdot j =n} a_k(b_j)\right) n^{-s}.$$ 
\end{fed} 

\begin{rem}
Let us suppose that $D$ is a $B(X)$-valued Dirichlet series, $E$ is an $X$-valued Dirichlet series, and the abscissa of convergence is finite for both series. Then there exists a half-plane $\mathbb{C}_\sigma$ in which both series converge, the first in $B(X)$ and the second in $X$.  Following what was done in the previous sections, the product could be defined in that half-plane as $D(s)(E(s))$. Let us see that this definition does not differ from Definition \ref{Def prod dirichlet}. That is, there exists $\sigma>0$ such that if $s$ is a complex number with real part greater than $\sigma$, then $D\cdot E(s) = D(s) (E(s))$ (and then in all the half-plane of convergence). For this, let $\sigma_0>0$ such that both series converge absolutely in $\mathbb{C}_{\sigma_0}$. Then,
$$D(s) (E(s)) = D(s) \left(\sum\limits_j b_j j^{-s}\right) = \sum\limits_j D(s)(b_j) j^{-s} = \sum\limits_j \left(\sum\limits_k a_k(b_j) k^{-s} \right) j^{-s},$$

It is sufficient to see that the series converge absolutely in $X$ to be able to reorder the sum. Since $\sum\limits_k \Vert a_k \Vert_{B(X)} k^{-\sigma} < \infty$ and $\sum\limits_j \Vert b_j\Vert_X j^{-\sigma} < \infty$ if $\sigma>\sigma_0$, then 

$$\sum\limits_n \Vert \sum\limits_{k j =n} a_k(b_j) \Vert_X n^{-\sigma} \leq \sum\limits_n \sum\limits_{kj=n} \Vert a_k\Vert_{B(X)}\Vert b_j\Vert_X n^{-\sigma} = \sum\limits_k \Vert a_k\Vert_{B(X)} k^{-\sigma} \sum\limits_{j} \Vert b_j \Vert_X j^{-\sigma} <\infty.$$

Then, if $s\in \mathbb{C}_{\sigma_0}$ we can rewrite:
$$D\cdot E(s) = \sum\limits_n \left(\sum\limits_{k\cdot j =n} a_k(b_j)\right) n^{-s} = \sum\limits_j \left(\sum\limits_k a_k(b_j) k^{-s} \right) j^{-s} = D(s) (E(s)).$$
\end{rem}
Formal power and Dirichlet series are connected by meaning of the Bohr transform $\mathcal{B}$. Given $n\in \mathbb{N}$, there exists a unique $\alpha\in \mathbb{N}_0^{(\mathbb{N})}$ such that $n = \mathfrak{p}^\alpha$, being $\mathfrak{p}$ the sequence of prime numbers and $\mathfrak{p}^\alpha = \mathfrak{p}_1^{\alpha_1}\cdot\mathfrak{p}_2^{\alpha_2} \cdots$. Then \begin{equation}\label{Def Bohr}
\mathcal{B}\left(\sum\limits_{\alpha\in \mathbb{N}_0^{(\mathbb{N})}} a_\alpha z^\alpha\right) = \sum\limits_{n=1}^\infty a_n n^{-s},
\end{equation}
where $a_n = a_\alpha$ if $n=\mathfrak{p}^\alpha$, defines a bijection between the formal power and Dirichlet series. In this way, if $f:\mathbb{T}^\infty\to X$ is formally associated with a power series $\sum_{\alpha\in \mathbb{N}_0^{(\mathbb{N})}} a_\alpha w^\alpha$, then by $\mathcal{B}(f)$ we refer to the Dirichlet series given by $\mathcal{B}(\sum_{\alpha\in \mathbb{N}_0^{(\mathbb{N})}}a_\alpha w^\alpha)$. Similarly, and without making any distinction, we will denote by $\mathcal{B}(F)$, for a holomorphic function $F:\mathbb{D}^\infty_2\to X$, the Dirichlet series obtained by applying $\mathcal{B}$ to the formal power series associated with $F$.

Now we  turn our attention to Hardy spaces of Dirichlet series with coefficients in a Banach space (see \cite{carando2014bohr,defantperez_2018} and \cite[24.2]{DirSer}). 
\begin{fed}
    Let $1\leq p \leq\infty$, the Hardy space $\mathcal{H}_p(X)$ is defined by
    $$\mathcal{H}_p(X) = \left\{D= \sum a_n n^{-s} : \mathcal{B}^{-1} (D)\sim f \in H_p(\mathbb{T}^\infty,X) \right\},$$
    which naturally is a Banach space with the norm $\Vert D \Vert_{\mathcal{H}_p(X)} = \Vert \mathcal{B}^{-1}(D) \Vert_{H_p(\mathbb{T}^\infty,X)}$.
\end{fed}

By the definition, the Dirichlet polynomials are dense, that is Dirichlet series with only finitely many non-zero coefficients. We also consider a second class of Hardy spaces of Dirichlet series (see \cite{defantperez_2018}).
\begin{fed}
    Let $1\leq p \leq \infty$, the Hardy space $\mathcal{H}_p^+(X)$ is defined by those Dirichlet series $D= \sum a_n n^{-s}$ such that $$D_\varepsilon = \sum\limits_n \frac{a_n}{n^\varepsilon} n^{-s} \in \mathcal{H}_p(X),$$  for all $ \varepsilon>0$ and $$\Vert D \Vert_{\mathcal{H}_p^+(X)} = \sup\limits_{\varepsilon>0} \Vert D_\varepsilon \Vert_{\mathcal{H}_p(X)} < \infty.$$
\end{fed}

Defant and Pérez prove in \cite{defantperez_2018} that these spaces are Banach spaces and isometrically contain $\mathcal{H}_p(X)$, for each $1\leq p \leq \infty$. Furthermore, in the same article, is proved that, for $1\leq p \leq \infty$, the Bohr transform




$$\mathcal{B}: H_p(\mathbb{D}^\infty_2,X) \to \mathcal{H}_p^+(X),$$
is an isometric isomorphism. In fact, for $p=\infty$ they actually show that $\mathcal{H}_\infty(X)$ is isomorphic to $H_\infty(B_{c_0},X)$, which in turn is isometrically isomorphic to $H_\infty(\mathbb{D}^\infty_2,X)$ as mentioned in Remark \ref{B_c0 a D_2}.

Then, we have that for all $1 \le p \le \infty$,  $\mathcal{H}_p(X)$ and $\mathcal{H}_p^+(X)$ are isometrically isomorphic if and only if $X$ has the ARNP. Furthermore, from Theorem \ref{iso con sot} we have that $\mathcal{H}_\infty^+(B(X))$ is isometrically isomorphic to $H_\infty^{sot}(\mathbb{T}^\infty,B(X))$ as long as $X$ is a separable Banach space with the ARNP. Our goal now  is to study the multipliers in the spaces $\mathcal{H}_p^+(X)$ and $\mathcal{H}_p(X)$.

\begin{teo}\label{Teo mult ser dir hol}
Let $X$ be a Banach space, $1\leq p < \infty$ and $D:\mathbb{C}_\sigma \to B(X)$ a function defined in some half-plane. Then the operator $M_D : \mathcal{H}^+_p(X) \to \mathcal{H}^+_p(X)$ given by $$M_D(E)(s) = D(s)(E(s))$$
is well-defined and bounded if and only if $D\in \mathcal{H}_\infty^+(B(X))$. Moreover, in this case $$\Vert D  \Vert_{\mathcal{H}_\infty^+(B(X))} = \Vert M_D \Vert.$$
\end{teo}

\bdem

Assume that $M_D : \mathcal{H}^+_p(X) \to \mathcal{H}^+_p(X)$ is well defined. Given $x \in X$ we have that $$M_D(x) = D(s)(x) \in \mathcal{H}^+_p(X),$$ then $D(s)(x) = \sum\limits_n a_n(x) n^{-s}$. We want to show  that $a_n : X \to X$, defined by $a_n(x)$ the $n-$th coefficient of $D(s)(x)$, is a bounded operator and that in some half-plane the equality $D(s) = \sum\limits_n a_n n^{-s}$ holds. If $\sigma$ is large enough, then for every $n$ and $x\in X$ we have

$$a_n(x) = \lim\limits_{R\to \infty}\frac{1}{2R}\int\limits_{-R}^R D(\sigma+it)(x) n^{\sigma+it}\mathrm{d}t.$$
The linearity of the integral and $D( \sigma+it)$ give the linearity of $a_n$. Since $$\Vert a_n(x) \Vert_X \leq \Vert D(s) (x)\Vert_{\mathcal{H}^+_p(X)} \leq \Vert M_D(x)\Vert_{\mathcal{H}^+_p(X)} \leq \Vert M_D \Vert \Vert x\Vert_{\mathcal{H}^+_p(X)} = \Vert M_D \Vert \Vert x \Vert_X,$$
we get that $a_n \in B(X)$ and $\Vert a_n \Vert_{B(X)}\leq \Vert M_D \Vert$. Then we can deduce that for $\sigma>1$ the sum $\sum\limits_n \frac{\Vert a_n \Vert_{B(X)}}{n^\sigma} < \infty$, and hence $\sum\limits_n  a_n n^{-s}$ converges absolutely in $B(X)$ on the half-plane $\mathbb{C}_1$. Then, if $x\in X$ and $s\in \mathbb{C}_1$ we have that
$$\left(\sum\limits_n  a_n n^{-s}\right)(x) = \sum\limits_n a_n(x) n^{-s} = D(s)(x),$$

hence  $D$ admits a Dirichlet series representation. It remains to show that $D\in \mathcal{H}_\infty^+(B(X))$.

Given a Dirichlet series $E$ we consider $E_N$ the series whose $n$-th coefficient coincides with that of $E$ if in the factorization of $n$ only the first $N$ primes appear and zero otherwise. Then for each $x\in X$ we have that $$D_N(s)(x) = \sum\limits_n a_n^N(x) n^{-s} \in \mathcal{H}_p^{+,(N)}(X) := \{ E_N : E\in \mathcal{H}^+_p(X)\},$$

and $\Vert D_N \Vert_{\mathcal{H}_p^{+,(N)}(X)} \leq \Vert D \Vert_{\mathcal{H}^+_p(X)}$. 
Then, the function 
$$F_N(z)(x) = \sum\limits_{\alpha\in \mathbb{N}_0^{N}} a_{\mathfrak{p}^\alpha}(x) z^\alpha \text{ belongs to\ }   
 H_p(\mathbb{D}^N,X).$$ Even more, since $\Vert a_n \Vert_{B(X)} \leq \Vert M_D\Vert$ for all $n\in \mathbb{N}$, we have that $\sum\limits_{\alpha\in \mathbb{N}_0^N} a_{\mathfrak{p}^\alpha} z^\alpha$ converges absolutely in $\mathbb{D}^N$. Then $F_N:\mathbb{D}^N \to B(X)$, given by $F_N(z) = \sum\limits_{\alpha\in \mathbb{N}_0^N} a_{\mathfrak{p}^\alpha} z^\alpha$, is holomorphic.

Let $G\in H_p(\mathbb{D}^N,X)$ and $E \in \mathcal{H}_p^{+,(N)}(X)$ the image by the Bohr transform of $G$, then $D\cdot E \in \mathcal{H}_p(X)$ and  $D_N\cdot E =(D \cdot E)_N\in \mathcal{H}_p^{+,(N)}(X)$. From looking at the coefficients we deduce that $F_N\cdot G \in H_p(\mathbb{D}^N,X)$. Therefore $F_N$ is a multiplier of $H_p(\mathbb{D}^N,X)$, and as in the infinite dimensional case we conclude that $F_N\in H_\infty(\mathbb{D}^N,B(X))$, that is, $D_N\in H_\infty^+(B(X))$.

If we consider  the inclusion $\iota : \mathcal{H}_p^{+,(N)}(X) \to \mathcal{H}^+_p(X)$ and the projection  $\Pi : \mathcal{H}^+_p(X) \to \mathcal{H}_p^{+,(N)}(X)$, we get $M_{D_N} = \Pi \circ M_D \circ \iota$, so $\Vert M_{D_N} \Vert \leq \Vert M_D\Vert$. Also $M_{F_N} = \mathcal{B}^{-1} \circ M_{D_N} \circ \mathcal{B}$, where $\mathcal{B}$ is the vector-valued Bohr transform. Then
$$\Vert F_N\Vert_{H_\infty(\mathbb{D}^N,B(X))} = \Vert M_{F_N} \Vert = \Vert M_{D_N} \Vert \leq \Vert M_D \Vert,$$  
hence there exists a function $F\in H_\infty(B_{c_0},B(X))  = H_\infty(\mathbb{D}^\infty_2,B(X))$ such that $C_\alpha(F) = a_{\mathfrak{p}^\alpha}$ for all $\alpha\in \mathbb{N}_0^{(N)}$ and $$\Vert F \Vert_{H_\infty(\mathbb{D}_2^\infty,B(X))} = \sup\limits_N \Vert F_N\Vert_{H_\infty(\mathbb{D}^N,B(X)) }\leq  \Vert M_D \Vert.$$

We have shown that $D\in \mathcal{H}_\infty^+(B(X))$ y $\Vert D \Vert_{\mathcal{H}_\infty^+(B(X))} \leq \Vert M_D\Vert.$

Suppose now that $D\in \mathcal{H}_\infty^+(B(X))$ and let $F = \mathcal{B}(D) \in H_\infty(\mathbb{D}_2^\infty,B(X))$. Then we can define the operator $M_F : H_p(\mathbb{D}_2^\infty,X) \to H_p(\mathbb{D}_2^\infty,X)$. Let $E\in \mathcal{H}^+_p(X)$ and $G= \mathcal{B}(E) \in H_p(\mathbb{D}_2^\infty,X)$, then $F\cdot G \in H_p(\mathbb{D}_2^\infty,X)$ and we have seen that the coefficients coincide with those of $D\cdot E$ (remember that $D\cdot E$ can be represented as a Dirichlet series in the half-plane of absolute convergence of both series). Applying again the Bohr transform we get that $D\cdot E \in \mathcal{H}^+_p(X)$. Then $M_D : \mathcal{H}^+_p(X) \to \mathcal{H}^+_p(X)$ is well defined and $$\Vert D \Vert_{\mathcal{H}^+_\infty(B(X))} = \Vert F \Vert_{H_\infty(\mathbb{D}_2^\infty,B(X))} = \Vert M_F \Vert  =\Vert M_D \Vert. $$
\edem

To study the multipliers of $\mathcal{H}_p(X)$, if $X$ is a separable Banach space, it is natural, by the definition of the spaces, to consider those formal Dirichlet series $D = \sum a_n n^{-s},$ with $a_n\in B(X)$, for which 
$$\sum a_{\mathfrak{p}^\alpha} z^\alpha \sim f \in  H_\infty^{sot}(\mathbb{T}^\infty,B(X)).$$
Note that every function $f$ in $H_\infty^{sot}(\mathbb{T}^\infty,B(X))$ defines, for each $\alpha\in \mathbb{N}_0^{(\mathbb{N})}$, a bounded operator $a_{\mathfrak{p}^\alpha} \in B(X)$, being $a_{\mathfrak{p}^\alpha}(x)$ the $\alpha-$th Fourier coefficient of $f(\cdot)(x) \in H_p(\mathbb{T}^\infty,X)$, and 
$$\Vert a_{\mathfrak{p}^\alpha}\Vert_{B(X)} \leq \Vert f \Vert_{H_\infty^{sot}(\mathbb{T}^\infty,B(X))}.$$
Therefore, the Dirichlet series $D = \sum a_n n^{-s}$ is well defined and for every $s\in \mathbb{C}_1$ then we have $D(s) \in B(X)$. We define the Banach space 
$$\mathcal{H}_\infty^{sot}(B(X)):=\{ D=\sum a_n n^{-s} : \sum a_{\mathfrak{p}^\alpha} z^\alpha \sim f \in  H_\infty^{sot}(\mathbb{T}^\infty,B(X))\},$$
with the norm 
$$\Vert D\Vert_{\mathcal{H}_\infty^{sot}(B(X))} = \Vert f \Vert_{H_\infty^{sot}(\mathbb{T}^\infty,B(X))}.$$

\begin{teo}\label{Teo mult ser dir toro}
    Let $X$ be a separable Banach space, $1\leq p < \infty$ and $D:\mathbb{C}_\sigma \to B(X)$ a function defined in some half-plane. Then the operator $M_D : \mathcal{H}_p(X) \to \mathcal{H}_p(X)$ given by $$M_D(E)(s) = D(s)(E(s))$$ is well-defined and bounded if and only if $D\in \mathcal{H}_\infty^{sot}(B(X))$. Moreover, in this case $$\Vert D  \Vert_{\mathcal{H}_\infty^{sot}(B(X))} = \Vert M_D \Vert.$$
\end{teo}
\bdem
Take $D=\sum a_n n^{-s}\in \mathcal{H}_\infty^{sot}(B(X))$ and let $f\in H_\infty^{sot}(\mathbb{T}^\infty,B(X))$ such that $f\sim \sum a_{\mathfrak{p}^\alpha} z^\alpha$, then $$M_f : H_p(\mathbb{T}^\infty,X) \to H_p(\mathbb{T}^\infty,X),$$
is well defined and $\Vert M_f \Vert = \Vert f \Vert_{H_\infty^{sot}(\mathbb{T}^\infty,B(X))}$. Given $E = \sum b_n n^{-s} \in \mathcal{H}_p(X)$ and $g \sim \sum b_{\mathfrak{p}^\alpha} z^{\alpha}$ the associated function by means of the Bohr transform. Then $f\cdot g \in H_p(\mathbb{T}^\infty,X)$ and using the density of the polynomials, we get that $$f\cdot g \sim \sum\limits_\alpha \left(  \sum\limits_{\beta+\gamma = \alpha} a_{\mathfrak{p^\beta}}(b_{\mathfrak{p}^\gamma})\right)z^\alpha.$$
Therefore, through the Bohr transform, we have $$D\cdot E =\sum\limits_{n} (\sum\limits_{k\cdot j =n} a_k(b_j)) n^{-s} \in \mathcal{H}_p(X),$$
and even more $\Vert D\cdot E \Vert_{\mathcal{H}_p(X)} = \Vert f\cdot g \Vert_{H_p(\mathbb{T}^\infty,X)}.$
We conclude that $M_D$ is well defined and $$\Vert M_D \Vert  = \Vert M_f \Vert = \Vert f\Vert_{H_\infty^{sot}(\mathbb{T}^\infty, B(X))} = \Vert D \Vert_{\mathcal{H}_\infty^{sot}(B(X))}.$$

Suppose now that $X$ is a separable Banach space and $D : \mathbb{C}_\sigma \to B(X)$ defines a multiplication operator. Then, the arguments of Theorem \ref{Teo mult ser dir hol} show that $D (s)$ admits a Dirichlet series representation $\sum\limits a_n n^{-s}$ on $\mathbb{C}_1$, and that for every $x\in X$
$$M_D(x) = \sum\limits a_n(x)n^{-s} \in \mathcal{H}_p(X).$$
Then $$\sum a_{\mathfrak{p}^\alpha}(x) z^{\alpha} \sim f_x \in H_p(\mathbb{T}^\infty,X),$$
and the linearity of $M_D$ implies that for every $x,y\in X$ and $\lambda,\beta\in \mathbb{C}$ the equality $f_{\lambda x+\beta y} = \lambda f_x + \beta f_y$ almost everywhere in $\mathbb{T}^\infty$. The same arguments given in Theorem \ref{iso con sot} prove that there exists $f\in H_\infty^{sot}(\mathbb{T}^\infty,B(X))$ such that $f(\cdot)(x) = f_x$ for every $x\in X$. Therefore $f\sim  \sum a_{\mathfrak{p}^\alpha} z^{\alpha}$ and hence $D\in \mathcal{H}_\infty^{sot}(B(X))$. From the first part we have that 
$$\Vert M_D \Vert  = \Vert D \Vert_{\mathcal{H}_\infty^{sot}(B(X))}.$$
\edem

From Theorems \ref{Teo mult ser dir hol} and \ref{Teo mult ser dir toro} we have the next remark
\begin{rem}

\begin{enumerate}
    \item[a)] If $X$ is a Banach space with the ARNP, then the space $\mathcal{H}_\infty^+(B(X))$ is the space of multipliers of $\mathcal{H}_p(X) = \mathcal{H}_p^+(X)$.
    \item[b)] If $X$ is separable but not has the ARNP, then the multipliers of $\mathcal{H}_p(X)$ and $ \mathcal{H}_p^+(X)$ differ, being $\mathcal{H}_\infty^{sot}(B(X))$ and $\mathcal{H}_\infty^+(B(X))$ respectively. 
    \item[c)]If $X$ is a separable Banach space with the ARNP, then $\mathcal{H}_\infty^+(B(X))$ and $\mathcal{H}_\infty^{sot}(B(X))$ coincide and are the multipliers of $\mathcal{H}_p(X)$. 
    \end{enumerate}
\end{rem}

Following Definition \ref{Def prod dirichlet}, we can consider as a multiplier of $\mathcal{H}_p(X)$ (or $\mathcal{H}_p^+(X)$) a formal Dirichlet series $D=\sum a_n n^{-s}$, with $a_n\in B(X)$. However, if $D$ is such a multiplier, then $M_D(x) = \sum a_n(x) n^{-s}$ belongs to $\mathcal{H}_p(X)$ (or $\mathcal{H}_p^+(X)$) and by the continuity of the coefficients, $\Vert a_n \Vert_{B(X)}$ is uniformly bounded. Then $D$ converges on $\mathbb{C}_1$ and the case is characterized in Theorems \ref{Teo mult ser dir hol} and \ref{Teo mult ser dir toro}. As a consequence, we have the following. 
\begin{pro}
    Let $X$ a Banach space, $D$ a formal Dirichlet series with coefficients in $B(X)$ and $1\leq p <\infty$ then the following statements hold.
    \begin{itemize}
        \item[i)] If $D$ is a multiplier on $\mathcal{H}_p^+(X)$ then $D\in \mathcal{H}_\infty^+(B(X))$.
        \item[ii)] If $X$ is separable and $D$ is a multiplier on $\mathcal{H}_p(X)$ then $D\in \mathcal{H}_\infty^{sot}(B(X))$.
    \end{itemize}
\end{pro}

Studying multipliers defined by formal Dirichlet series is, without doubt, simpler than studying those given by formal power series in the spaces $H_\infty(\mathbb{T}^\infty,X)$ and $H_\infty(\mathbb{D}^\infty_2,X)$. This is due to the fact that working with Dirichlet series involves one-variable complex analysis, where convergence at a single point ensures convergence and holomorphy throughout an entire half-plane. Because of the connection between formal Dirichlet series and formal power series through the Bohr transform we get our last result. 

\begin{pro}\label{prop-last}
    Let $X$ a Banach space, $\varphi=\sum a_\alpha z^{\alpha}$ a formal power series and $1\leq p<\infty$ then the following statements hold
    \begin{itemize}
        \item[i)] If $\varphi$ is a multiplier on $H_p(\mathbb{D}_2^\infty,X)$, then $\varphi\sim F$ for some $F\in H_\infty(\mathbb{D}_2^\infty,B(X))$.
        \item[ii)] If $X$ has the ARNP and $\varphi$ is a multiplier on $H_p(\mathbb{T}^\infty,X)$, then $\varphi\sim F$ for some $F\in H_\infty(\mathbb{D}_2^\infty,B(X))$.
        \item[iii)] If $X$ is separable and $\varphi$ is a multiplier of  $H_p(\mathbb{T}^\infty,X)$, then $\varphi\sim f$ for some $f\in H_\infty^{sot}(\mathbb{T}^\infty,B(X))$. 
    \end{itemize}
\end{pro}
\bibliographystyle{abbrvurl}
\bibliography{refs}

\noindent
J.~Antezana\\
Departament de Matemàtiques i Informàtica, Universitat de Barcelona. Barcelona, Spain\\
Centre de Recerca Matemàtica, Barcelona, Spain.

D.~Carando, T.~Fern\'andez Vidal\\
Departamento de Matem\'{a}tica,
Facultad de Cs. Exactas y Naturales, Universidad de Buenos Aires and IMAS-CONICET. C.A.B.A., Argentina, dcarando@dm.uba.ar, tfvidal@dm.uba.ar\\ 

M.~Scotti\\
Departamento de Matem\'{a}tica,
Facultad de Cs. Exactas y Naturales, Universidad de Buenos Aires. C.A.B.A., Argentina, mscotti@dm.uba.ar\\

\end{document}